\documentclass[12pt]{amsart}

\usepackage{amssymb}
\usepackage{latexsym}
\usepackage{verbatim}
\usepackage{fullpage}

\input {cyracc.def}
\font\tencyr=wncyr10 %scaled \magstephalf
\font\tencyi=wncyi10 %scaled \magstephalf
\font\tencysc=wncysc10 %scaled \magstephalf
\def\rus{\tencyr\cyracc}
\def\rusi{\tencyi\cyracc}
\def\rusc{\tencysc\cyracc}

\newcommand{\re}[1]{\textrm  (\ref{#1})}

\renewenvironment{proof}
{\noindent {\sl Proof.}\quad }{\hfill
$\square$ \vskip1.1ex\noindent }

\newenvironment{proof*}
{\noindent {\sl Proof.}\quad }{\hfill
$\square$}

\renewcommand{\theequation}{\thesection .\arabic{equation}}
\renewcommand{\thesubsubsection}{\theequation .\arabic{subsubsection}}

\catcode`\@=\active
\catcode`\@=11
\def\@eqnnum{\hbox to
.01pt{}\rlap{\hskip-\displaywidth(\mathbf{\theequation})}}
\catcode`\@=12
\newenvironment{s}[1]
{ \vskip1.2ex \refstepcounter{equation}
\noindent {\bf \theequation\enspace #1.} \begin{sl}}{\end{sl}
\vskip1.1ex\noindent }

\newenvironment{rem}[1]
{ \vskip1.2ex \refstepcounter{equation}
\noindent {\bf \theequation\enspace {#1}.} }{ \vskip1.1ex\noindent }

\newenvironment{subs}[1]
{\vskip1.2ex \refstepcounter{equation}
\noindent {\bf (\theequation)\quad #1.} }{\quad}

\newcommand {\ah}{{\frak a}}
\newcommand {\be}{{\frak b}}
\newcommand {\ce}{{\frak c}}

\newcommand {\g}{{\frak g}}

\newcommand {\el}{{\frak l}}
\newcommand {\me}{{\frak m}}
\newcommand {\n}{{\frak n}}

\newcommand {\p}{{\frak p}}

\newcommand {\te}{{\frak t}}
\newcommand {\ut}{{\frak u}}

\newcommand {\sln}{{\frak sl}_n}
\newcommand {\sltn}{{\frak sl}_{2n}}
\newcommand {\slno}{{\frak sl}_{n+1}}

\newcommand {\spn}{{\frak sp}_{2n}}
\newcommand {\sono}{{\frak so}_{2n+1}}
\newcommand {\sone}{{\frak so}_{2n}}

\newcommand {\gA}{{\goth A}}
\newcommand {\gC}{{\goth C}}
\newcommand {\gH}{{\goth H}}

\newcommand {\ap}{\alpha}

\newcommand {\lb}{\lambda}
\newcommand {\vp}{\varphi}

\newcommand {\HW}{\widehat W}
\newcommand {\HV}{\widehat V}
\newcommand {\HP}{\widehat\Pi}
\newcommand {\HD}{\widehat\Delta}

\newcommand {\hot}{{\mathrm{ht\,}}}

\newcommand {\rk}{{\mathrm{rk\,}}}

\newcommand {\srk}{{\mathrm{srk\,}}}

\newcommand {\tri}{{\frak sl}_2}
\newcommand {\GR}[2]{{\textrm{{\bf #1}}}_{#2}}

\newcommand {\ov}{\overline}
\newcommand {\un}{\underline}

\newcommand {\Ab}{{\frak Ab}}
\newcommand {\AD}{{\frak Ad}}

\newcommand {\beq}{\begin{equation}}
\newcommand {\eeq}{\end{equation}}

\newcommand{\Shid}{\mathrm{Shi}(\Delta)}

\newcommand{\curle}{\preccurlyeq}

\newfam\Bbbfam\newfam\eufam\newfam\eusfam%
\font\Bbbfont=msbm10 scaled 1200%
\font\olala=msam10 scaled 1200%
\font\Bbbsmallfont=msbm8%
\textfont\Bbbfam=\Bbbfont\scriptfont\Bbbfam=\Bbbsmallfont%%
\font\euzw=eufm10 scaled 1200%
\font\euac=eufm7 scaled 1200%
\font\euacc=eufm7 scaled 1000%
\textfont\eufam=\euzw\scriptfont\eufam=\euac%
\scriptscriptfont\eufam=\euacc%
\font\euszw=eusm10 scaled 1200%
\font\eusac=eusm7 scaled 1200%
\font\eusacc=eusm7 scaled 1000%
\textfont\eusfam=\euszw\scriptfont\eusfam=\eusac%
\scriptscriptfont\eusfam=\eusacc%
\def\frak{\fam\eufam}%
\def\goth{\fam\eusfam}%
\def\Bbb{\fam\Bbbfam}%

\def\varnothing{\hbox {\Bbbfont\char'077}}
\def\square{\hbox {\olala\char"03}}

\def\cyeq{\hbox {\olala\char'064}}

\begin{document}
\setlength{\parskip}{2pt plus 4pt minus 0pt}
\hfill {\scriptsize May 19, 2004}
\vskip1ex
\vskip1ex

\title[Normalizers of ad-nilpotent ideals]{Normalizers of ad-nilpotent ideals}
%\medskip \\
\author{\sc Dmitri I. Panyushev}
\thanks{This research was supported in part by 
%R.F.B.I. Grant no. 02--01--01041 and 
CRDF Grant no. RM1-2543-MO-03 and SFB/TR12 at the Ruhr-Universit\"at Bochum}
\maketitle
\begin{center}
{\footnotesize
{\it Independent University of Moscow,
Bol'shoi Vlasevskii per. 11 \\
119002 Moscow, \quad Russia \\ e-mail}: {\tt panyush@mccme.ru }\\
}
\end{center}

\noindent
Let $\g$ be a complex simple Lie algebra.
Fix a Borel subalgebra $\be$ and a Cartan subalgebra $\te\subset\be$.
The nilpotent radical of $\be$ is denoted by $\ut$.
The corresponding set of positive (resp. simple) roots is $\Delta^+$
(resp. $\Pi$).
\\[.4ex]
An ideal of $\be$ is called
{\sf ad}-{\it nilpotent\/}, if it is contained in $[\be,\be]$.
The theory of {\sf ad}-nilpotent ideals has attracted much recent
attention in the work of Kostant, Cellini-Papi, Sommers, and others.
The goal of this paper is to study the normalizer of an
{\sf ad}-nilpotent ideal. We obtain several general descriptions of the
normalizer, and present a number of more explicit results
for $\g=\sln$ or $\spn$.

Let $\ce$ be an {\sf ad}-nilpotent ideal. Being a
$\te$-stable subspace of $\ut$, it is a sum of root spaces.
The sum of the corresponding roots is an integral weight,
denoted $|\ce|$.
%%which is said to be the {\it weight\/} of $\ce$.
We show that $|\ce|$ is a dominant weight and that
the normalizer of $\ce$, $\n_\g(\ce)$, is completely determined by
$|\ce|$.
Since $\n_\g(\ce)$ contains $\be$, it suffices to realize
which root spaces $\g_{-\ap}$ ($\ap\in\Pi$) are contained in
$\n_\g(\ce)$. We prove that
$\g_{-\ap}\subset \n_\g(\ce)$ if and only if $(|\ce|,\ap)=0$.

Another type of descriptions is based on a relationship between the
{\sf ad}-nilpotent ideals and certain elements of the affine Weyl group
$\HW$. Let ${\AD}=\AD(\g)$ denote the set of all {\sf ad}-nilpotent ideals
of $\be$.
By \cite{CP1}, to each $\ce\in\AD$ one associates an element
of $\HW$, which we denote by $w_{min,\ce}$.
 An {\sf ad}-nilpotent ideal is called {\it strictly positive\/},
if it is contained in $[\ut,\ut]$. The set of strictly positive ideals
is denoted by $\AD_0$.
By \cite{eric}, to each $\ce\in\AD_0$ one associates an element of $\HW$,
which we denote by $w_{max,\ce}$. The group $\HW$ acts linearly on
a vector space $\HV$, containing affine root system, and
we prove that
\[
\g_{-\ap}\subset \n_\g(\ce) \ \Longleftrightarrow \
w_{min,\ce}(\ap) \text{ is an affine simple root } .
\]
If $\ce\in\AD_0$, then both $w_{min,\ce}$ and $w_{max,\ce}$ are defined, and we have
\[
w_{min,\ce}(\ap) \text{ is an affine simple root }
\ \Longleftrightarrow \
w_{max,\ce}(\ap) \text{ is an affine simple root } ,
\]
see Sections~\ref{prelim} and \ref{sect-norm} for details.
It is worth noting that for $\ce\in\AD_0$ the elements
$w_{min,\ce}$ and $w_{max,\ce}$ can be different and simple roots
$w_{min,\ce}(\ap)$, $w_{max,\ce}(\ap)$ can also be different.
\\[.5ex]
Our geometric characterization of the normalizer is connected with Shi's bijection
between the {\sf ad}-nilpotent ideals and the dominant regions of the Shi arrangement.
Given an ideal $\ce\in\AD$, let $R_\ce$ denote the corresponding region. We show
that $\g_{-\ap}\subset \n_\g(\ce)$ if and only if the wall of the dominant Weyl
chamber orthogonal to $\ap$ is also a wall of $R_\ce$.

An interesting general problem is to study the partition of $\AD$ into the subsets
$\AD\{\p\}$ parametrized by the standard parabolic subalgebras.
Here $\AD\{\p\}:=\{\ce\in\AD \mid \n_\g(\ce)=\p\}$.
Obviously, the nilpotent radical of $\p$ is the unique maximal element
of $\AD\{\p\}$. But there can be several minimal elements.
It seems that the most interesting subset is $\AD\{\be\}$.
We give a description of it in the spirit of Cellini-Papi.
According to \cite{CP2}, there is a bijection between $\AD$ and
the points of the coroot lattice $Q^\vee$ lying in the simplex
$D_{min}=\{x\in V\mid (\ap,x)\ge -1 \ \forall \ap\in\Pi \ \& \ (\theta,x)\le 2\}$.
Here $V=\underset{\ap\in\Pi}{\oplus} {\Bbb R}\ap$, $\theta$ is the highest
root, and  $(\ ,\ )$ is a $W$-invariant inner product on $V$.
Our result says that the normalizer of the ideal corresponding 
to $x\in D_{min}\cap Q^\vee$ is $\be$ if and only if
$(x,\ap)\ne 0$ for all $\ap\in\Pi$ and $(x,\theta)\ne 1$.
%%integral point of a certain simplex, and we explicitly indicate
%%those points that correspond to $\AD\{\be\}$.
This description allows us to interpret the number $\#\AD\{\be\}$ as
the coefficient of $x$ in a certain Laurent series. 
This series depends only on the coefficients of $\theta$ in the basis of 
the simple roots. Namely, if
$\theta=\sum_{i=1}^n c_i\ap_i$, $c_0=1$, and $f=\#\{j\mid c_j=1\}$, then
\[
  \#\AD\{\be\}=\frac{1}{f}[x]\prod_{i=0}^n (\frac{x^{-c_i}}{1-x^{c_i}}-1) \ .
\]
A similar characterization is obtained for $\AD_0\{\be\}:=\AD\{\be\}\cap\AD_0$.
Using these results we explicitly compute the numbers \ 
$\#\AD\{\be\}$ and $\#\AD_0\{\be\}$ for the classical Lie algebras. 

In case of $\g=\slno$ and $\spn$, our results are more precise.
We explicitly describe the set $\AD\{\p\}$ for any  $\p$.
For $\g=\slno$, it follows from that
description that $\#\AD\{\p\}$ depends only on the difference $s=n-\srk\p$.
Namely, it is the {\it $s$-th Motzkin number}.
For $\g=\spn$, the quantity $\#\AD\{\p\}$ depends only on the number of short simple
roots, say $s$, that are not in the Levi subalgebra of $\p$.
Namely, it is
the {\it number of directed animals of size\/} $s+1$.
In these two cases, it is also shown that $\AD\{\p\}$ has always a unique minimal
element. The interest of two series, $\slno$ and $\spn$, is revealed via the fact
that one can tie together the notion of duality for {\sf ad}-nilpotent ideals
\cite{duality},
the  minimax ideals \cite{mima}, and the ideals whose normalizer equals $\be$.
That is, we prove that $\n_\g(\ce)=\be$ if and only if the dual ideal, $\ce^*$,
is minimax. Also, it turns out that considering the normalizers of {\sf ad}-nilpotent
ideals provides a natural framework for demonstrating various identities
related to Catalan, Motzkin, Riordan, and other numbers.

{\small {\bf Acknowledgements.} A part of this paper was written during
my stay at the Ruhr-Universit\"at Bochum.
I would like to thank Prof. A.T.\,Huckleberry for invitation and
hospitality. Thanks are also due to G.\,R\"ohrle for interesting
discussions on {\sf ad}-nilpotent ideals.}

                %%%%
%%%%%%%%%%%%%%%%
%%%%%%%%%%%%%%%%   Notation
%%%%%%%%%%%%%%%%
                %%%%

\section{Notation and other preliminaries}
\label{prelim}

\noindent
\begin{subs}{Main notation}
\end{subs}
$\Delta$ is the root system of $(\g,\te)$ and
$W$ is the usual Weyl group. For $\mu\in\Delta$, $\g_\mu$ is the
corresponding root space in $\g$.

$\Delta^+$  is the set of positive roots, $\theta$ is the highest root in $\Delta^+$,
and $\rho=\frac{1}{2}\sum_{\mu\in\Delta^+}\mu$.

$\Pi=\{\ap_1,\dots,\ap_p\}$ is the set of simple roots in $\Delta^+$
and 
$\vp_i$ is the fundamental weight corresponding to $\ap_i$.
%
%%$e_1,e_2,\dots,e_p$ \ are the exponents and $h$ is the Coxeter
%%number of $W$.
 \\
 We set $V:=\te_{\Bbb R}=\oplus_{i=1}^p{\Bbb R}\ap_i$ and denote by
$(\ ,\ )$ a $W$-invariant inner product on $V$. As usual,
$\mu^\vee=2\mu/(\mu,\mu)$ is the coroot
for $\mu\in \Delta$.

${\gC}=\{x\in V\mid (x,\ap)>0 \ \ \forall \ap\in\Pi\}$
\ is the (open) fundamental Weyl chamber.

${\gA}=\{x\in V\mid (x,\ap)>0 \ \  \forall \ap\in\Pi \ \ \& \
(x,\theta)<1\}$ \ is the fundamental alcove.

%%$Q=\oplus _{i=1}^p {\Bbb Z}\ap_i \subset V$ is the root lattice

$Q^+=\{\sum_{i=1}^p n_i\ap_i \mid n_i=0,1,2,\dots \}$
and $Q^\vee=\oplus _{i=1}^p {\Bbb Z}\ap_i^\vee\subset V$
is the coroot lattice.
\\
For $\mu,\gamma\in\Delta^+$, write $\mu\curle\gamma$, if
$\gamma-\mu\in Q^+$.
%%The notation $\mu\prec\gamma$ means that
%%$\mu\curle\gamma$ and $\gamma\ne\mu$.
We regard $\Delta^+$ as poset under `$\cyeq$'.
\\[.5ex]
Letting $\widehat V=V\oplus {\Bbb R}\delta\oplus {\Bbb R}\lb$, we extend
the inner product $(\ ,\ )$ on $\widehat V$ so that $(\delta,V)=(\lb,V)=
(\delta,\delta)=
(\lb,\lb)=0$ and $(\delta,\lb)=1$.

$\widehat\Delta=\{\Delta+k\delta \mid k\in {\Bbb Z}\}$ is the set of affine
real roots and $\widehat W$ is the  affine Weyl group.
\\
Then $\widehat\Delta^+= \Delta^+ \cup \{ \Delta +k\delta \mid k\ge 1\}$ is
the set of positive
affine roots and $\widehat \Pi=\Pi\cup\{\ap_0\}$ is the corresponding set
of affine simple roots,
where $\ap_0=\delta-\theta$.
%%, where $\theta$ is the highest root in $\Delta^+$.
The inner product $(\ ,\ )$ on $\widehat V$ is
$\widehat W$-invariant. The notation $\beta>0$ (resp. $\beta <0$)
is a shorthand for $\beta\in\HD^+$ (resp. $\beta\in -\HD^+$).
\\
For $\ap_i$ ($0\le i\le p$), we let $s_i$ denote the corresponding simple
reflection in $\widehat W$.
If the index of $\ap\in\widehat\Pi$ is not specified, then we merely write
$s_\ap$. %% for the corresponding reflection.
The length function on $\widehat W$ with respect
to  $s_0,s_1,\dots,s_p$ is denoted by $\ell$.
For any $w\in\widehat W$, we set
\[
   N(w)=\{\mu\in\widehat\Delta^+ \mid w(\mu) \in -\widehat \Delta^+ \} .
\]
It is standard that $\#N(w)=\ell(w)$ and $N(w)$ is {\it bi-convex\/}.
The latter means
that both $N(w)$ and $\HD^+\setminus N(w)$ are subsets of $\HD^+$
that are closed under addition.
Furthermore, the assignment $w\mapsto N(w)$ sets up
a bijection between the elements of $\HW$ and the finite bi-convex subsets
of $\HD^+$.

\begin{subs}{Ideals and antichains}               \label{ideals}
\end{subs}
Our $\be$ is the Borel subalgebra of $\g$ corresponding
to $\Delta^+$ and $\ut=[\be,\be]$.
If $\ce$ is an {\sf ad}-nilpotent ideal of $\be$, then
$\ce=\underset{\ap\in I_\ce}{\oplus}\g_\ap$
for some $I_\ce\subset \Delta^+$. The set of roots $I_\ce$ ($\ce\in\AD$)
arising in this way is an {\it upper ideal\/} of (the poset)
$\Delta^+$. This means that $I_\ce$ satisfies the following property:\\
\centerline{
if $\gamma\in I_\ce,\mu\in\Delta^+$, and $\gamma+\mu\in\Delta$, then
$\gamma+\mu\in I_\ce$.}
In view of the obvious bijection between the {\sf ad}-nilpotent ideals and
the upper ideals of $\Delta^+$, we will often identify them.
A root $\gamma\in I_\ce$ is called
a {\it generator of\/} $I_\ce$ or $\ce$, if $\gamma-\ap\not\in I_\ce$ for any
$\ap\in\Pi$.
In other words, $\gamma$ is a minimal element of $I_\ce$ with respect to ``$\curle$".
We write $\Gamma(I_\ce)$ or $\Gamma(\ce)$ for the set of generators.
It is easily seen that $\Gamma(I_\ce)$ is an {\it antichain\/}
of $\Delta^+$, i.e., $\gamma_i\not\curle\gamma_j$ for any pair
$(\gamma_i,\gamma_j)$ in $\Gamma(I_\ce)$.
Conversely, if $\Gamma\subset \Delta^+$ is an antichain,
then the upper ideal
$I\langle\Gamma\rangle:=\{\mu\in \Delta^+\mid \mu\succcurlyeq\gamma_i
\text{ for some } \gamma_i\in \Gamma \}$
has $\Gamma$ as the set of generators.

                %%%%
%%%%%%%%%%%%%%%%
%%%%%%%%%%%%%%%%   Section 2  
%%%%%%%%%%%%%%%%
                %%%%

\section{The weight and normalizer of an {\sf ad}-nilpotent ideal}  
\label{sect-norm}
\setcounter{equation}{0}

\noindent
A parabolic subalgebra of $\g$ is called {\it standard\/}, if it contains
$\be$. If $\Pi'\subset \Pi$, then $\p(\Pi')$ stands for the standard parabolic
subalgebra which is generated by $\be$ and the spaces $\g_{-\ap}$, $\ap\in\Pi'$.
For instance, $\p(\varnothing)=\be$ and $\p(\Pi)=\g$. 
The maximal (proper) parabolic subalgebra $\p(\Pi\setminus \{\ap_i\})$ is
also denoted by $\p_{\langle i\rangle}$,
and we write $\p(\ap_i)$ in place of $\p(\{\ap_i\})$.
The standard Levi subalgebra of $\p(\Pi')$, denoted $\el(\Pi')$,
is generated by $\te$ and the subspaces $\g_\ap$, $\ap\in\pm\Pi'$.
%%The number  equals
Write $\srk\p(\Pi')$ for $\rk [\el(\Pi'),\el(\Pi')]$,
the semisimple rank of $\p(\Pi')$.
We have $\srk\p(\Pi')=\#\Pi'$.
\\[.6ex]
Let $\n_\g(\ce)$ be the normalizer in $\g$ of an {\sf ad}-nilpotent ideal $\ce$.
It is a standard parabolic subalgebra of $\g$.
Therefore, to describe $\n_\g(\ce)$ explicitly, one has to only
realize when $\g_{-\ap}$ is contained in $\n_\g(\ce)$ for an $\ap\in \Pi$.
A description of $\n_\g(\ce)$ in terms of $\Gamma(\ce)$ is given in
\cite[Theorem\,3.2]{pr}:

\begin{s}{Theorem}  \label{p-r}
$\g_{-\ap}\not\subset\n_\g(\ce)$ is and only if $\gamma-\ap\in\Delta^+\cup
\{0\}$ for some $\gamma\in\Gamma(\ce)$.
\end{s}%
The aim of this section is to give some other descriptions of 
$\n_\g(\ce)$ associated with the combinatorial theory of
{\sf ad}-nilpotent ideals.

\noindent
Recall some basic results concerning a connection between
the {\sf ad}-nilpotent ideals and certain elements in the affine Weyl group.
Given $\ce\in \AD$ with the corresponding upper ideal $I_\ce\subset\Delta^+$, there is
a unique element $w=w_{min,\ce}\in\HW$ satisfying the
following properties (see \cite{CP1}):
\begin{enumerate}
\item[($\Diamond$)] \ For $\gamma\in \Delta^+$, we have $\gamma\in I_\ce$ if and only
if $w(\delta-\gamma) < 0$;
\item[\sf (dom)] \ $w(\ap)>0$ for all $\ap\in\Pi$;
\item[\sf (min)] \ if $\ap\in\HP$ and  
$w^{-1}(\ap)=k\delta+\mu$ for some $\mu\in \Delta$, then $k\ge -1$.
\end{enumerate}
This element $w$ is said to be the {\it minimal element of\/} $\ce$.
The elements of $\HW$ satisfying  property {\sf (dom)} are called {\it dominant}.
The elements of $\HW$ satisfying the last two properties are called {\it minimal}.
The minimal element of $\ce$ can also be characterized as the unique element
of $\HW$ satisfying properties ($\Diamond$), {\sf (dom)}, 
and having the minimal possible length.
This explains the term.
\\[.6ex]
An {\sf ad}-nilpotent ideal $\ce$ is called {\it strictly positive\/}, if it is 
contained in $[\ut,\ut]$ 
(i.e., $I_\ce$ contains no simple roots). 
The set of strictly positive ideals is denoted by $\AD_0$.
If $\ce\in\AD_0$, then there is a unique element 
$w=w_{max,\ce}\in\HW$ satisfying properties ($\Diamond$) and {\sf (dom)},
as above, and also the property
%%following properties (see \cite{eric}):
\begin{enumerate}
%%\item \ $w(\ap)>0$ for all $\ap\in\Pi$;
\item[\sf (max)] \  if $\ap\in\HP$ and  
$w^{-1}(\ap)=k\delta+\mu$ for some $\mu\in \Delta$, then $k\le 1$,
\end{enumerate}
(see \cite{eric}).
This element is said to be the {\it maximal element of\/} $\ce$.
The elements of $\HW$ satisfying properties {\sf (dom)} and {\sf (max)}
are called {\it maximal}.
The maximal element of a strictly positive ideal 
can also be characterized as the unique element
of $\HW$ satisfying properties ($\Diamond$), {\sf (dom)}, and having the maximal
possible length. This explains the term.
\\
Usually, we have $w_{min,\ce}\ne w_{max,\ce}$. The case of coincidence of these
two elements is studied in \cite{mima}. The respective ideals are called 
{\it minimax\/}.
\\[.6ex]
For any finite subset $M\subset \widehat\Delta^+$, we set\
$\vert M\vert:=\sum_{\gamma\in M}\gamma$. 
%%Clearly, it is an element of $Q^+$. 
If $\ce\in \AD$ and $I_\ce$ is the corresponding
upper ideal, then put $|\ce|:=|I_\ce|$.
We say that $|\ce|$ is the {\it weight\/} of the {\sf ad}-nilpotent ideal $\ce$.
Our first aim is to look at the weights occurring in this way.
The following result is due to Kostant \cite[Theorem\,7]{top}. 
For the sake of completeness, we give a proof, which demonstrates the
role of minimal elements.

\begin{s}{Proposition} \label{uniqueness}
Suppose $\ce_1,\ce_2\in\AD$ and $|\ce_1|=|\ce_2|$. Then $\ce_1=\ce_2$.
\end{s}\begin{proof}
Let $I_1,I_2$ be the corresponding upper ideals.
Assume $I_1\ne I_2$. Then both sets $I_1\setminus I_2$ and $I_2\setminus I_1$
are non-empty and we have $|I_1\setminus I_2|=|I_2\setminus I_1|$.
Let us rewrite this equality in the form:
\begin{equation}   \label{rewrite}
  \sum_{\gamma\in I_1\setminus I_2}(\delta-\gamma)-c\delta=
   \sum_{\gamma\in I_2\setminus I_1}(\delta-\gamma) \ ,
\end{equation}
where $c=\#I_1- \#I_2$.
Without loss of generality, we may assume that $\dim \ce_1\ge \dim\ce_2$,
i.e., $c\ge 0$. 
Let $w_1\in \HW$ be the minimal element of $\ce_1$.
Applying $w_1$ to Eq.~\re{rewrite} and using property ($\Diamond$), we see that
$w_1(\text{LHS})$ (resp. $w_1(\text{RHS})$) is a sum of negative (resp. positive)
roots. A contradiction!
\end{proof}%%
For $\ce\in\AD$, we set $\ce^1=\ce$ and $\ce^k=[\ce^{k-1},\ce]$ for $k\ge 2$. 
Then $\ce^m=0$ for $m\gg 0$.

\begin{s}{Theorem}  \label{ves}
Let $\ce$ be an {\sf ad}-nilpotent ideal of $\be$ and $\ap\in\Pi$. Then
\begin{itemize}
\item[\sf (i)] \ $(|\ce|,\ap)\ge 0$ ;
\item[\sf (ii)] \ $\g_{-\ap}\subset\n_\g(\ce) \ \ 
\Leftrightarrow \ (|\ce|,\ap)=0$;
\item[\sf (iii)] \  $(|\ce|,\ap)=0 \ \Leftrightarrow \ (\sum_{k\ge 1}|\ce^k|,\ap)=0$.
\end{itemize}
\end{s}\begin{proof}
(i),\,(ii) For $\ap\in\Pi$, 
let $\tri(\ap)$ be the simple three-dimensional subalgebra of $\g$ generated
by $\g_\ap$ and $\g_{-\ap}$. Let $\{x_\ap,h_\ap,y_{-\ap}\}$ be a basis for 
$\tri(\ap)$, where $x_\ap\in \g_\ap$, $y_{-\ap}\in \g_{-\ap}$, and
$h_\ap=[x_\ap,y_{-\ap}]$. 
Obviously, $\ce$ is a $\langle x_\ap,h_\ap\rangle$-module. 
Since $\ce$ is a subspace of an $\tri(\ap)$-module, we conclude that
$(|\ce|,\ap)\ge 0$. This proves part (i). 
Furthermore, $(|\ce|,\ap)=0$ if and only if $\ce$ is an $\tri(\ap)$-module,
i.e., $y_{-\ap}\in \n_\g(\ce)$.

(iii) Since $\ce^k$ is an {\sf ad}-nilpotent ideal for all $k\ge 1$, we have
$(|\ce^k|,\ap) \ge 0$ by part (i).
This gives the implication ``$\Leftarrow$". On the other hand,
if $\g_{-\ap}\subset\n_\g(\ce)$, then $\g_{-\ap}\subset\n_\g(\ce^k)$
as well, and one may apply part (ii) to $\ce^k$.
\end{proof}%%
Thus, the weight of any ideal is dominant, different ideals have different weights,
and the normalizer of an {\sf ad}-nilpotent
ideal is completely determined by its weight.

\begin{rem}{Remarks}
1. If $\ce_1$ and $\ce_2$ are two {\sf ad}-nilpotent ideals, then 
$\ce_1\cap\ce_2$ and $\ce_1+\ce_2$ are {\sf ad}-nilpotent ideals as well.
Also, $|\ce_1+\ce_2|+|\ce_1\cap\ce_2|=|\ce_1|+|\ce_2|$.
Clearly,
\[
   \n_\g(\ce_1+\ce_2)\supset \n_\g(\ce_1)\cap \n_\g(\ce_2)  
\ \text{ and } \ 
    \n_\g(\ce_1\cap\ce_2)\supset \n_\g(\ce_1)\cap \n_\g(\ce_2) \ .
\]
But both these containments can be strict even if $\n_\g(\ce_1)=\n_\g(\ce_2)$,
see Example~\ref{primery}(2) below.

2. It is an interesting open problem to characterize abstractly  the set
of weights $\{ |\ce| \mid \ce\in\AD\}$.
For instance, if $\g=\GR{G}{2}$, then this set is equal to 
$\{0,\,\vp_2,\, 3\vp_1,\,4\vp_1,\,3\vp_1+\vp_2,\, 5\vp_1,\,3\vp_2,\,2\vp_1+2\vp_2\}$.
\end{rem}%
We wish also to obtain a description of $\n_\g(\ce)$ in terms of $w_{min,\ce}$ 
(and $w_{max,\ce}$, if $\ce\in\AD_0$). To this end, consider
\[
   \hat\rho=\rho+ \frac{(\theta,\theta)}{2}(1+(\rho,\theta^\vee))\lb \in \HV \ .
\]
Since $(\rho, \ap^\vee)=1$ for any $\ap\in\Pi$,
it readily follows that $\hat\rho$ is the unique element of $V\oplus{\Bbb R}\lb$
having the property that $(\hat\rho,\ap^\vee)=1$ for all $\ap\in\HP$.

\begin{s}{Proposition}  \label{razn}
For any $w\in \HW$, we have\/ $\hat\rho-w^{-1}(\hat\rho)=|N(w)|$.
\end{s}\begin{proof*}
We argue by induction on $\ell(w)$. If $w=s_\ap$, $\ap\in\HP$,
then $N(s_\ap)=\{\ap\}$ and the claim follows from the
definition of $\hat\rho$.
Suppose $\ell(w)>1$ and $w=s_\ap \tilde w$, where $\ell(w)=\ell(\tilde w)+1$.
Then $N(w)=N(\tilde w)\cup \{\tilde w^{-1}(\ap)\}$. Assume that the claim holds for 
$\tilde w$.
Then $\hat\rho-w^{-1}(\hat\rho)=\hat\rho-\tilde w^{-1}(\hat\rho)+\tilde w^{-1}(\hat\rho)
-w^{-1}(\hat\rho)=|N(\tilde w)|+\tilde w^{-1}(\hat\rho-s_\ap\hat\rho)=
|N(\tilde w)|+\tilde w^{-1}(\ap)=|N(w)|$.
\end{proof*}%
\begin{s}{Lemma}   \label{maxim1}
Let $w\in\HW$ 
be a dominant element. Then 
\begin{itemize}
\item[\sf (i)] \ $(|N(w)|,\ap)\le 0$ for any $\ap\in\Pi$;
\item[\sf (ii)] \ 
$(|N(w)|,\ap)=0$ if and only if $w(\ap)\in\HP$;
\end{itemize}
\end{s}\begin{proof}
By Proposition~\ref{razn}, $\hat\rho-w^{-1}(\hat\rho)=|N(w)|$.
It follows that $(|N(w)|,\ap^\vee)=1-(\hat\rho,w(\ap)^\vee)$.
By property ({\sf dom}), we have $w(\ap)$ is positive.
This yields all the assertions.
\end{proof}%
The following is our main result for minimal elements.

\begin{s}{Theorem}   \label{norm-min}
Let $\ce$ be an arbitrary\/ {\sf ad}-nilpotent ideal of $\be$ and $\ap\in\Pi$.
Then $\g_{-\ap}\subset\n_\g(\ce)$ if and only if $w_{min,\ce}(\ap)\in\HP$.
\end{s}\begin{proof}
By \cite[Section\,2]{CP1}, we have 
$N(w_{min,\ce})=\displaystyle\bigcup_{k\ge 1}(k\delta-I_{\ce^k}$).
Hence, $(|N(w_{min,\ce})|,\ap)=
-(\sum_{k\ge 1} |\ce^k|,\ap)$. Therefore,
combining Theorem~\ref{ves} and Lemma~\ref{maxim1}, we obtain:\\
\centerline{
$\g_{-\ap}\subset\n_\g(\ce)$ \ if and only if \ 
$(|N(w_{min,\ce})|,\ap)=0$ \ if and only if \
$w_{min,\ce}(\ap)\in\HP$.}
\end{proof}%
Next, we show that, for a strictly positive ideal $\ce$, the similar claim holds
with $w_{max,\ce}$.

\begin{s}{Theorem}   \label{maxim2}
Suppose $\ce\in\AD_0$ and $\ap\in\Pi$. Then
%\centerline{
$w_{max,\ce}(\ap)\in\HP$ if and only if $w_{min,\ce}(\ap)\in\HP$.
%}
\end{s}\begin{proof}
We have already proved that the two conditions:
\begin{itemize}
\item \ $w_{min,\ce}(\ap)\in\HP$,
%%\item \ $(|N(w_{min})|,\ap)=0$,
\item \ $\ce$ is an $\tri(\ap)$-module
\end{itemize}
are equivalent.
Therefore, by Lemma~\ref{maxim1}, it suffices to prove that
$(|N(w_{max,\ce})|,\ap)=0$ if and only if $\ce$ is an $\tri(\ap)$-module.
A description of $N(w_{max,\ce})$ is due to Sommers \cite{eric},
see also \cite[2.11]{losh}. We state it in a form convenient for our purposes.
Let $\me$ be the (unique) $\te$-stable complement of $\ce$ in $\ut$.
Set $\me^1=\me$ and $\me^k=[\me^{k-1},\me]$ for $k\ge 2$.
Let $\tilde\ce^k$ be the $\te$-stable complement of $\me+\me^2+\ldots+\me^k$
in $\ut$.
Then $\ce=\tilde\ce\supset\tilde\ce^2\supset\ldots$ and $\tilde\ce^m=0$ for $m\gg 0$.
By \cite{eric}, we have $\ce^k\subset\tilde\ce^k$ for any $k$ and
$\displaystyle  N(w_{max,\ce})=\bigcup_{k\ge 1}(k\delta-I_{\tilde\ce^k})$.
Hence
\begin{equation}  \label{zvezdochka}
(|N(w_{max,\ce})|,\ap)=-\sum_k(|\tilde\ce^k|,\ap) \ .
\end{equation}
Notice that $\g_\ap\subset\me$.
Let $\me_\ap$ be the $\te$-stable complement of $\g_\ap$ in $\me$.
Then $\me_\ap\oplus\ce=\p(\ap)^{nil}$, the nilpotent radical of the minimal parabolic
subalgebra $\p(\ap)$. Since $\p(\ap)$ is an $\tri(\ap)$-module and $\ce$ is an
$\langle x_\ap,h_\ap\rangle$-module, $\me_\ap$ is an
$\langle y_\ap,h_\ap\rangle$-module. Furthermore, $\me_\ap$ is an
$\tri(\ap)$-module if and only if $\ce$ is. Next, we have
$\me+\me^2+\ldots+\me^k=\g_\ap\oplus (\me_\ap+\me^2+\ldots+\me^k)$.
Set $\me_\ap^{\langle k\rangle}:=(\me_\ap+\me^2+\ldots+\me^k)$ for $k\ge 2$.
Then $\me_\ap^{\langle k\rangle}\oplus\tilde\ce^k=\p(\ap)^{nil}$.
By induction, it is easily seen that each $\me_\ap^{\langle k\rangle}$ is
an $\langle y_\ap,h_\ap\rangle$-module. Therefore $\tilde\ce^k$ is an
$\langle x_\ap,h_\ap\rangle$-module and hence $(|\tilde\ce^k|, \ap)\ge 0$
for all $k\ge 1$.
If $(|N(w_{max,\ce})|,\ap)=0$, then it follows from Eq.~\re{zvezdochka} that
$(|\tilde\ce^k|, \ap)=0$ for all $k$. In particular, $\tilde\ce=\ce$ is an
$\tri(\ap)$-module.
Conversely, if $\ce$ is an $\tri(\ap)$-module, then
$\me_\ap$ is. This easily implies that
$\me_\ap^{\langle k\rangle}$ is an $\tri(\ap)$-module for all $k\ge 2$.
Hence all $\tilde\ce^k$,  $k\ge 1$, are $\tri(\ap)$-modules and we conclude from
Eq.~\re{zvezdochka} that $(|N(w_{max,\ce})|,\ap)=0$.
\end{proof}%
It is not however true in general that two simple roots
$w_{max,\ce}(\ap)\in\HP$ and $w_{min,\ce}(\ap)\in\HP$ are equal, see
Example~\ref{primery}(4) below.

\begin{s}{Corollary}  \label{norm-max}
For $\ce\in\AD_0$ and $\ap\in\Pi$, we have $\g_{-\ap}\subset\n_\g(\ce)$ if and only if\/
$(\sum_k|\tilde\ce^k|,\ap)=0$ if and only\/ if $w_{max,\ce}(\ap)\in\HP$.
\end{s}%
\vskip-1.2ex
\begin{rem}{Examples}   \label{primery}
Here we give some illustrations to previous results.
The numbering of the simple roots and fundamental weights is the same is in \cite{VO}.

{\sf (1)} Let $\ce$ be the {\sf ad}-nilpotent ideal for $\g={\frak sl}_5$ with
$\Gamma(\ce)=\{\ap_1+\ap_2,\,\ap_2+\ap_3+\ap_4\}$. It is an Abelian ideal, i.e.,
$\ce^2=0$. Here
$|\ce|=2\vp_1+2\vp_2+\vp_4$. Using either Theorem~\ref{p-r} or Theorem~\ref{ves},
we obtain $\n_\g(\ce)=\p(\ap_3)$. Using the algorithm given in 
\cite[6.1]{lp}, one finds
that $w_{min,\ce}=s_1s_4s_5s_0$.
Therefore the action of $w_{min,\ce}$ on $\Pi$ is given by
\[
     w_{min,\ce}:\left\{\begin{array}{ccl}
%%%\ap_0 & \mapsto & -\delta+\ap_2 \\
\ap_1 & \mapsto & \delta-\ap_1-\ap_2\\
\ap_2 & \mapsto & \ap_1+\ap_2+\ap_3\\
\ap_3 & \mapsto & \ap_4\\
\ap_4 & \mapsto & \delta-\ap_2-\ap_3-\ap_4
\end{array}\right.
\]
Then using Theorem~\ref{norm-min} we obtain again the same description of
$\n_\g(\ce)$.
In this case $\ce^2=0$, but $\tilde\ce^2=\g_{\theta}$ and $\tilde\ce^3=0$.
Therefore $\sum_k|\tilde\ce^k|=|\ce|+\theta=3\vp_1+2\vp_2+2\vp_4$, which yields
an illustration to Corollary~\ref{norm-max}.
Finally, one can compute that $w_{max,\ce}=s_2w_{min,\ce}$.

{\sf (2)}  Take $\g={\frak sl}_7$ and the ideals $\ce_1,\ce_2$ with
$\Gamma(\ce_1)=\{\ap_1{+}\ap_2,\ap_3{+}\ap_4{+}\ap_5,\ap_4{+}\ap_5{+}\ap_6\}$,
$\Gamma(\ce_2)=\{\ap_1{+}\ap_2{+}\ap_3,\ap_2{+}\ap_3{+}\ap_4,\ap_4{+}\ap_5,
\ap_5{+}\ap_6\}$. Then one easily computes that
$\n_\g(\ce_1)=\n_\g(\ce_2)=\be$, whereas 
$\n_\g(\ce_1\cap\ce_2)=\p(\ap_2)$ and $\n_\g(\ce_1+\ce_2)=\p(\ap_3)$.

{\sf (3)} $\g=\GR{F}{4}$. 
\\
Write $[n_1,\,n_2,\,n_3,\,n_4]$ for 
$\sum_i n_i\ap_i$.
Consider the ideal $\ce$ with 
$\Gamma(\ce)=\{[0,2,2,1],[2,2,1,0]\}$. 
Here $|\ce|=[16,28,20,10]=4\vp_1+2\vp_3$.
Next, 
$\ce^2=\tilde\ce^2=\g_{[2,4,3,1]}\oplus\g_{[2,4,3,2]}$ and $\tilde\ce^3=0$.
Therefore $\sum_k|\ce^k|=\sum_k|\tilde\ce^k|=[20,36,26,13]=4\vp_1+3\vp_3$.
Thus, $\n_\g(\ce)=\p(\{\ap_2,\ap_4\})$.
In this case, we have $w_{min,\ce}=w_{max,\ce}=s_0
s_4s_3s_2s_0s_4s_3s_1s_2s_3s_4s_0$.

{\sf (4)} $\g=\GR{G}{2}$. 
\\
Consider the Abelian ideal $\ce$ with $\Gamma(\ce)=\{2\ap_1+\ap_2\}$.
Here $w_{min,\ce}=s_1s_2s_0$ and $w_{max,\ce}=s_0s_2s_1s_2s_0$.
Therefore
\[
     w_{min,\ce}:\left\{\begin{array}{ccl}
%%%\ap_0 & \mapsto & -\delta+\ap_2 \\
\ap_1 & \mapsto & 2\ap_1+\ap_2\\
\ap_2 & \mapsto & \delta-3\ap_1-2\ap_2=\ap_0
\end{array}\right. , \quad
w_{max,\ce}:\left\{\begin{array}{ccl}
%%%\ap_0 & \mapsto & -\delta+\ap_2 \\
\ap_1 & \mapsto & \delta-\ap_1-\ap_2\\
\ap_2 & \mapsto & \ap_2
\end{array}\right. \ .
\]
Thus, we have $\n_\g(\ce)=\p(\ap_2)$, but $w_{min,\ce}(\ap_2)\ne w_{max,\ce}(\ap_2)$.
\end{rem}%
%

                %%%%
%%%%%%%%%%%%%%%%
%%%%%%%%%%%%%%%%   Section 3
%%%%%%%%%%%%%%%%
                %%%%

\section{A geometric description of the normalizer}  
\label{walls}
\setcounter{equation}{0}

\noindent
For $\mu\in\Delta^+$ and $k\in{\Bbb Z}$, define the hyperplane
$\gH_{\mu,k}$ in $V$ as $\{x\in V \mid (x,\mu)=k\}$.
The collection of all these hyperplanes is called the {\it affine arrangement\/}
in $V$.
The {\it regions\/} of an arrangement are the connected components of
the complement in $V$ of the union of all its hyperplanes. 
As is well-known (see e.g. \cite{hump}), the regions of the affine arrangement
are simplices ({\it alcoves}), and $\gA$ is one of them.
\\[.5ex]
Recall a bijection between the {\sf ad}-nilpotent ideals (or
upper ideals of $\Delta^+$) and the dominant regions of the Shi arrangement,
which is due to J.-Y.\,Shi \cite[Theorem\,1.4]{shi}.
%%%, see also \cite[Lemma\,3.1]{ath03}.
The {\it Shi arrangement\/}, $\mathrm{Shi}(\Delta)$,
is the sub-arrangement of the affine arrangement consisting of all hyperplanes
$\gH_{\mu,k}$ with
%%$\mu\in \Delta^+$ and
$k=0,1$.
Any region of $\mathrm{Shi}(\Delta)$ lying in $\gC$ is said to be
{\it dominant\/}.
The Shi bijection takes an {\sf ad}-nilpotent ideal $\ce$ with
corresponding upper ideal $I_\ce\subset \Delta^+$ to the dominant region
\begin{equation}   \label{bij-shi}
 R_\ce =\{ x\in \gC \mid  (x,\gamma)>1, \text{ if \ }\gamma\in I_\ce \quad \&\quad
           (x,\gamma)<1,  \text{ if \ }\gamma\not\in I_\ce
           \}\ .
\end{equation}
Our goal is to describe $\n_\g(\ce)$ in terms of $R_\ce$.
To this end, we use relations between the two actions of $\HW$:
the linear action on $\HV$ and the affine action on $V$.
We use `$\ast$' to denote the affine action:
$(w,x)\mapsto w{\ast} x$, $w\in\HW,\ x\in V$.
For any $\ap\in\HP$, let $H_\ap$ denote the corresponding wall of $\gA$.
That is, $H_\ap=
\left\{\begin{array}{cl} \gH_{\ap,0}, & \text{ if } \ \ap\in\Pi, \\
                         \gH_{\theta,1}, & \text{ if } \ \ap=\ap_0
\end{array}\right.$.
%%Recall that $\HW$ naturally acts on $V$ as a group of affine transformations.
The generator $s_\ap\in\HW$ acts on $V$ as (affine) reflection relative to
$H_\ap$.
Our next arguments will be based on comparing properties of these two actions.
%%The closure of the fundamental alcove $\gA$ is a fundamental domain for this action.
The following is Eq.~(1.1) in \cite{CP1}. Suppose $\mu\in\Delta^+$,
$k>0$, and $h\ge 0$. Then
\begin{equation}  \label{eq-cp1}
 \begin{array}{rlc}  w(k\delta-\mu)<0 & \text{ if and only if }\ \gH_{\mu,k}
   & \text{ separates } \ \gA \text{ and } \ w^{-1}{\ast}\gA \ , \\
      w(h\delta+\mu)<0 & \text{ if and only if }\ \gH_{\mu,-h}
   & \text{ separates } \ \gA \text{ and } \ w^{-1}{\ast}\gA \ .
  \end{array}
\end{equation}
It~follows from these equations
%%\cite[(1.1)]{CP1}
that $w\in \HW$ is dominant if and only if $w^{-1}{\ast}\gA\subset\gC$.
Another useful relation is
\begin{equation}  \label{another}
w{\ast}\gH_{\mu_1,k_1}=\gH_{\mu_2,k_2} \
\text{ if and only if } \ w(k_1\delta-\mu_1)=\pm(k_2\delta-\mu_2) \ ,
\end{equation}
where $\mu_i\in\Delta^+$ and $k_i\in\Bbb Z$.
It suffices to verify this only for the simple reflections,
the case of $w=s_i$ $(i>0)$ being trivial. Some calculations are only needed
for $w=s_0$.
\\[.5ex]
Notice that if $R$ is a dominant region, then its walls belong to the
set $\gH_{\gamma,1}$, $\gamma\in \Delta^+$, and $\gH_{\ap,0}$, $\ap\in\Pi$.

\begin{s}{Theorem}  \label{shi-geom}
Suppose $\ce\in\AD$ and $\ap\in\Pi$. Then
$\g_{-\ap}\subset \n_\g(\ce)$ if and only if\/ $\gH_{\ap,0}$ is a wall of $R_\ce$.
\end{s}\begin{proof}
For any $w\in \HW$, let ${\goth L}(w)$ denote the set of all hyperplanes
$\gH_{\gamma,k}$ separating $\gA$ and $w{\ast}\gA$.

Suppose $\g_{-\ap}\subset \n_\g(\ce)$. Then $w_{min,\ce}(\ap)=:\nu\in\HP$.
Then
%%$w^{-1}(\nu)=\ap$ and therefore
$N(s_{\nu}w_{min,\ce})=N(w_{min,\ce})\cup\{\ap\}$. 
This already means that $\tilde w:=s_\nu w_{min,\ce}$ is not 
dominant. Furthermore, by Theorem\,4.5 in \cite{hump}, we have
${\goth L}(\tilde w^{-1})={\goth L}(w_{min,\ce}^{-1})\cup 
\{w_{min,\ce}^{-1}{\ast}H_{\nu}\}$ and by Eq.~\re{another},
$w_{min,\ce}^{-1}{\ast}H_{\nu}=\gH_{\ap,0}$.
That is, the hyperplane $\gH_{\ap,0}$ separates the
alcoves $\tilde w^{-1}{\ast}\gA$ and $w_{min,\ce}^{-1}{\ast}\gA$.
Since $w_{min,\ce}^{-1}{\ast}\gA\subset R_\ce$ \cite{CP2}, we conclude that
$\gH_{\ap,0}$ is a wall of $R_\ce$. 

Conversely, suppose $\gH_{\ap,0}$ is a wall of $R_\ce$. This  means 
that there is a $w\in\HW$ such that $w^{-1}{\ast}\gA\subset R_\ce$ 
(hence $w$ is dominant!) and $\gH_{\ap,0}$
is a wall of the alcove $w^{-1}{\ast}\gA$. Equivalently, $w{\ast}\gH_{\ap,0}$
is a wall of $\gA$. Then $w{\ast}\gH_{\ap,0}=H_\nu$ for some $\nu\in\HP$ and
hence $w(\ap)=\pm\nu$, by Eq.~\re{another}. Since $w$ is dominant, we actually
have $w(\ap)=\nu$. Next, it follows from the dominance of $w$ that
$N(w)=\bigcup_{k\ge 1}(k\delta- I_k)$, where 
each $I_k$ is an upper ideal. Furthermore, in view of Eq.~\re{eq-cp1},
the condition
$w^{-1}{\ast}\gA\subset R_\ce$ precisely means that $\delta-\gamma\in N(w)$
if and only if $\gamma\in I_\ce$, i.e., $I_1=I_\ce$.
Since $w(\ap)\in\HP$, it follows from Lemma~\ref{maxim1} that 
$(|N(w)|,\ap)=0$ and hence $(|I_k|,\ap)=0$ for all $k$.
In particular, $(|I_1|,\ap)=(|\ce|,\ap)=0$, i.e.,
$\g_{-\ap}\subset \n_\g(\ce)$.
\end{proof}%
{\bf Remark.} For $\ce\in\AD_0$, the previous result and Theorem~\ref{maxim2} 
show that the alcoves $w_{min,\ce}^{-1}{\ast}\gA$
and $w_{max,\ce}^{-1}{\ast}\gA$ have the same walls of the form
$\gH_{\ap,0}$. Furthermore, the following theorem shows that
if $\widetilde\gA$ is an arbitrary alcove in $R_\ce$, where $\ce$ is not
necessarily in $\AD_0$, 
then any its wall of the form $\gH_{\ap,0}$ is also a wall of
$w_{min,\ce}^{-1}{\ast}\gA$. 
\begin{s}{Theorem}  \label{subtle}
Suppose $w\in\HW$ is dominant, and let $\ce$ be the first layer ideal
of $w$ (i.e., $\ce=\{\gamma\in\Delta^+\mid w(\delta-\gamma)<0\}$). 
If $w(\nu)\in\HP$ for $\nu\in\Pi$, then $\g_{-\nu}\subset \n_\g(\ce)$.
\end{s}\begin{proof}
Assume $\g_{-\nu}\not\subset \n_\g(\ce)$. Then there is a $\gamma\in I_\ce$ such that
either $\gamma-\nu\in\Delta^+\setminus I_\ce$ or $\gamma=\nu$. In the first case we have
%%This means that 
$w(\delta-\gamma)<0$ and $w(\delta-\gamma+\nu)>0$. This clearly implies that
$\hot w(\nu)\ge 2$, i.e., this root is not simple.
If $\gamma=\nu$, then $\kappa:=w(\delta-\gamma)<0$. Hence
$w(\nu)=\delta-\kappa$. This root also cannot be simple.
\end{proof}%
%
%\\[.7ex]
Theorem~\ref{shi-geom} says that $\n_\g(\ce)\supset \p(\ap)$ if and only if 
$\gH_{\ap,0}$ is a wall of $R_\ce$. This can also be restated in the following form.
Consider the restricted arrangement
\[
\Shid_{\ap}=\{H\cap\gH_{\ap,0} \mid H\in\Shid\setminus \{\gH_{\ap,0}\cup\gH_{\ap,1}\}
\,\} \ .
\]
Let us say that the
region of $\Shid_{\ap}$ is dominant, if it belongs to $\bar\gC\cap \gH_{\ap,0}$.
Hence the ideals $\ce$ with $\n_\g(\ce)\supset \p(\ap)$ are in a bijection
with the dominant regions of $\Shid_{\ap}$.
Notice also that $\gH_{\gamma,1}\cap\gH_{\ap,0}=\gH_{\tilde\gamma,1}\cap\gH_{\ap,0}$
if and only if $\gamma-\tilde\gamma=k\ap$ for some $k\in \Bbb Z$.
This means the hyperplanes of the restricted arrangement that dissect
$\bar\gC\cap \gH_{\ap,0}$ are in a bijection with $\el(\ap)$-submodules of
$\p(\ap)^{nil}$.
(Here $\el(\ap)=\tri(\ap)+\te$ is the standard Levi subalgebra
of $\p(\ap)$.)
On the other hand, any {\sf ad}-nilpotent ideal whose normalizer
contains $\p(\ap)$ lies in $\p(\ap)^{nil}$ and is a sum of $\el(\ap)$-modules.
In the general case, the condition that a $\te$-stable subspace of
$\ut=\be^{nil}$ is actually
$\be$-stable led us to the notion of an upper ideal of $\Delta^+$.
Accordingly, in this situation, the condition that an $\el(\ap)$-stable
subspace of $\p(\ap)^{nil}$ is actually
$\p(\ap)$-stable lead us to the notion of an upper ideal of the poset of
{\it $\el(\ap)$-modules in $\p(\ap)^{nil}$}. The latter can be defined as the
quotient $\Delta^+_\ap:=(\Delta^+\setminus\{\ap\})/{\sim}$, where
the equivalence
$\sim$ is defined as $\gamma\sim\tilde\gamma$ if and only if
$\gamma-\tilde\gamma=k\ap$ for some $k\in \Bbb Z$. It is easily seen that "$\curle$"
induces a well-defined partial order in $\Delta^+_\ap$.
Thus, we obtain a ``restricted'' version of the Shi correspondence:
\\[.8ex]
\hbox to \textwidth{\refstepcounter{equation} {\bf (\theequation)} \hfil
\parbox{424pt}{\it There is a bijection between the upper ideals of $(\Delta^+_\ap,\,\curle)$ and
the dominant regions of the restricted arrangement $\Shid_\ap$.}
}
\\[.7ex]
Clearly, one can proceed further, and consider arbitrary parabolic
subalgebras (i.e., not necessarily minimal ones) and
the restricted Shi arrangement
determined by the respective face of the dominant Weyl chamber.
We leave it to the interested reader to 
give an accurate statement.
%%fill in the technical details.
It would be interesting to find a closed formula for the number of such 
dominant regions.

\vskip1ex
%\noindent
Let ${\frak Par}={\frak Par}(\g)$
denote the set of all standard parabolic subalgebras of $\g$.
We have a natural mapping $\psi: \AD\to {\frak Par}$, which takes an {\sf ad}-nilpotent
ideal to its normalizer. It is interesting to study the fibres of $\psi$.
Write $\AD\{\p\}$ for $\psi^{-1}(\p)$, the set of all ideals whose
normalizer equals $\p$.
Whenever we wish to make the dependence on $\g$
explicit, we write $\AD(\g)\{\p\}$.
Each $\AD\{\p\}$, as well as the whole of $\AD$, is regarded as poset under
the usual containment of subspaces of $\ut$.
The following is obvious.

\begin{s}{Lemma}   \label{fibre-psi}
The unique maximal element of $\AD\{\p\}$ is 
$\p^{nil}$, the nilpotent radical of $\p$. In particular,
$\psi$ is onto.
\end{s}%
It is not however true that $\AD\{\p\}$ always has a unique
minimal element.

{\bf Example.} Take $\g={\frak so}_8$ and $\p=\be$. Then
$\AD({\frak so}_8)\{\be\}$ has three minimal elements (ideals).
One of them has the generators $\ap_1{+}\ap_2, \ap_2{+}\ap_3{+}\ap_4$.
The other two correspond to the cyclic permutations of $\{1,3,4\}$.
\\[.5ex]
Below, we show that if $\g=\sln$ or $\spn$, then $\AD\{\p\}$
has a unique minimal element for any $\p\in {\frak Par}$.

Another easy observation is connected with the maximal parabolic subalgebras.

\begin{s}{Lemma}   \label{fibre-max}
The poset $\AD\{\p_{\langle i\rangle}\}$ is a chain and 
$\#\AD\{\p_{\langle i\rangle}\}=n_i$, where $\theta=\sum_i n_i\ap_i$.
%$\psi^{-1}(\p_{\langle i\rangle})=n_i$.
\end{s}\begin{proof}
For any $\mu\in Q$, let $[\mu:\ap_i]$ denote the coefficient of $\ap_i$ in
the expansion of $\mu$ via the basis $\Pi$.
Set $I(\ap_i)_j=\{\mu\in \Delta^+ \mid [\mu:\ap_i] \ge j\}$.
It is an upper ideal, and it is easily seen that the corresponding
{\sf ad}-nilpotent ideals with $j=1,\ldots,n_i$
exhaust the fibre $\AD\{\p_{\langle i\rangle}\}$.
%%$\psi^{-1}(\p_{\langle i\rangle})$.
\end{proof}%
In the rest of the section, we give a geometric description of the set $\AD\{\be\}$.
Recall that $\widehat W$ is isomorphic to a semi-direct product
of $W$ and $Q^\vee$ \cite{hump}.
Given $w\in\widehat W$, there is a unique factorization
\begin{equation}  \label{affine}
w=v{\cdot}t_{r}\ ,
\end{equation}
where $v\in W$ and  $t_{r}$ is the translation
corresponding to $r\in Q^\vee$. Then $w^{-1}=v^{-1}{\cdot}t_{-v(r)}$.
In terms of this factorization for $w$, the linear action of $\HW$ on 
$V\oplus {\Bbb R}\delta\subset\HV$
is given by 
$w^{-1}(x)=v^{-1}(x)+(x, v(r))\delta$ \ for any
$x\in V\oplus {\Bbb R}\delta$. In particular,
\begin{equation} \label{simroots} \left\{\begin{array}{l}
  w^{-1}(\ap_i)=v^{-1}(\ap_i)+(\ap_i, v(r))\delta, \quad i\ge 1 ,\\
  w^{-1}(\ap_0)=-v^{-1}(\theta)+(1-(\theta, v(r)))\delta \ .
\end{array}\right.
\end{equation}
Given $\ce\in\AD$, consider $w_{min,\ce}$ and the corresponding factorization
\re{affine} for it. (To simplify notation, we do not endow the components
$v$ and $r$ with subscripts.) Form the element $z_\ce:=v(r)\in Q^\vee\subset V$.
The following fundamental result is due to Cellini and Papi \cite{CP2}.

\begin{s}{Theorem}   \label{main-cp2}
\begin{itemize}
\item[\sf (i)] \ $z_\ce\in D_{min}:=
\{x\in V \mid (x,\ap)\ge -1 \ \ \forall\ap\in\Pi \ \ \& \
\ (x,\theta)\le 2\}$;
\item[\sf (ii)] \ The mapping $\AD \to D_{min}\cap Q^\vee$,
$\ce\mapsto z_\ce$, is a bijection.
\end{itemize}
\end{s}%
Our next description of $\AD\{\be\}$ says which points of $D_{min}\cap Q^\vee$
correspond to the ideals whose normalizer is equal to $\be$.

\begin{s}{Theorem}   \label{ad-b}  For $\ce\in\AD$, we have
\\
{\sf (i)} \ $\n_\g(\ce)=\be$ if and only if
%%\left\{\begin{array}{ll}
$(z_\ce, \ap)\ne 0 \ \ \forall\,\ap\in\Pi$ and
                          $(z_\ce,\theta)\ne 1$.
%%%\end{array}\right.$
In other words, there is a one-to-one correspondence
\[
   \AD\{\be\} \overset{1:1}{\longleftrightarrow} \{x\in D_{min}\cap Q^\vee \mid
   x\not\in H_\ap \ \ \forall\,\ap\in\HP \} .
%%(z_\ce, \ap)\ne 0 \ \ \forall\,\ap\in\Pi \ \ \& \ \ (z_\ce,\theta)\ne 1 \} .
\] 
{\sf (ii)} \ The semisimple rank of $\n_\g(\ce)$ is equal to the number of hyperplanes $H_\ap$ ($\ap\in\HP$)
to which $z_\ce$ belongs.
\end{s}\begin{proof}
(i) \  By Theorem~\ref{norm-min}, we have
\[
  \ce\in\AD\{\be\} \ \Longleftrightarrow \ w_{min,\ce}(\ap)\not\in\HP\ \
\forall\,\ap\in\Pi
\ \Longleftrightarrow \ w_{min,\ce}^{-1}(\HP)\cap\Pi=\varnothing \ .
\]
Actually, we will prove a bit more precise statement that

for $\ap\in\Pi$, \ $w_{min,\ce}^{-1}(\ap)\in\Pi$ if and only if $(z_\ce,\ap)=0$,

for $\ap=\ap_0$, $w_{min,\ce}^{-1}(\ap_0)\in\Pi$ if and only if $(z_\ce,\theta)=1$,
\\[.5ex]
which implies the assertion. Indeed, if $\ap\in\Pi$ and $w_{min,\ce}^{-1}(\ap)\in\Pi$,
then it follows from the first row in Eq.~\re{simroots}, with $z_\ce=v(r)$,
that $v^{-1}(\ap)\in\Pi$ and
$(z_\ce,\ap)=0$. Conversely, if $(z_\ce,\ap)=0$, then $w_{min,\ce}^{-1}(\ap)=
v^{-1}(\ap)=:\gamma\in\Delta$. Since $w_{min,\ce}$ is dominant, $\gamma$ must be positive.
Assuming $\gamma\not\in\Pi$ and hence $\gamma=\gamma_1+\gamma_2$ ($\gamma_i\in\Delta^+$),
we obtain $w_{min,\ce}(\gamma_1)+w_{min,\ce}(\gamma_2)=\ap$. Here both summands in the
left-hand side are positive roots, which contradicts the simplicity of $\ap$.
Hence $\gamma$ must be a simple root.

The argument for $\ap_0$ is similar, taking into account the second row in
Eq.~\re{simroots}.

(ii)  This follows from the proof of part (i).
If $z_\ce\in H_\ap$ ($\ap\in\HP$), then the corresponding root of 
$\Pi$ occurring in the standard Levi subalgebra of $\n_\g(\ce)$ is
$w_{min,\ce}(\ap)$. 
%%(But it is not yet clear how to describe this root entirely
%%in terms of $z_\ce$.)
\end{proof}%
It would be interesting to
%%%In principle, one should
be able to extract all information on the normalizer of
$\ce$ directly from the indication of $z_\ce$, without appealing to
$w_{min,\ce}$.
So that Theorem~\ref{ad-b} does not give a complete
answer to this problem.
\\[.6ex]
If $\ce\in\AD_0$, then similar results are valid for $w_{max,\ce}$ and the
coroot lattice points of another simplex.
For $w_{max,\ce}=v{\cdot}t_r$, we set $y_\ce=v(r)$.
The following fundamental result is due to Sommers \cite{eric}.

\begin{s}{Theorem}   \label{main-eric}
\begin{itemize}
\item[\sf (i)] \ $y_\ce\in D_{max}:=
\{x\in V \mid (x,\ap)\le 1 \ \ \forall\ap\in\Pi \ \ \& \
\ (x,\theta)\ge 0\}$;
\item[\sf (ii)] \ The mapping $\AD_0 \to D_{max}\cap Q^\vee$,
$\ce\mapsto y_\ce$, is a bijection.
\end{itemize}
\end{s}%
Our next description of $\AD_0\{\be\}$ says which points of $D_{max}\cap Q^\vee$
correspond to the ideals whose normalizer is equal to $\be$.
Since the proof is identical to the proof of Theorem~\ref{ad-b},
it is omitted.

\begin{s}{Theorem}   \label{ad0-b}
For $\ce\in\AD_0$, we have  \\
%%\begin{itemize} \item
{\sf (i)} \ $\n_\g(\ce)=\be$ if and only if
$(y_\ce, \ap)\ne 0 \ \ \forall\,\ap\in\Pi$ and
                          $(y_\ce,\theta)\ne 1$.
In other words, there is a one-to-one correspondence
\[
   \AD_0\{\be\} \overset{1:1}{\longleftrightarrow} \{x\in D_{max}\cap Q^\vee \mid
   x\not\in H_\ap \ \ \forall\,\ap\in\HP \} .
\]
{\sf (ii)} \ The semisimple rank of $\n_\g(\ce)$ is equal to the number
of hyperplanes $H_\ap$ ($\ap\in\HP$)
to which $y_\ce$ belongs.
\end{s}%
%

                %%%%
%%%%%%%%%%%%%%%%
%%%%%%%%%%%%%%%%   Section 4
%%%%%%%%%%%%%%%%
                %%%%

\section{On ideals whose normalizer is equal to $\be$}
\label{counting}
\setcounter{equation}{0}

\noindent 
In this section, we present a practical method for counting the number of
all and strictly positive {\sf ad}-nilpotent ideals,
respectively, whose normalizer equals $\be$.
It turns out that for the classical series of simple Lie algebras we meet 
several famous
integer sequences: the Motzkin and Riordan numbers, the number of directed 
animals of size $n$, and trinomial coefficients. 
Our exposition is quite similar to that in Section~5 in \cite{mima}, where 
an analogous problem was considered for minimax ideals.

%%%%%%%%%%%%%%%%%%%
\begin{comment}
By Theorem~\ref{ad-b}, we have
\[
\#\AD\{\be\}=\{ x\in Q^\vee \mid (x,\ap)\in\{-1,1,2,\ldots \} \ \ \forall \ap\in\Pi 
\quad \& \quad (x,\theta)\in \{2,0,-1,\ldots \} \}\ . 
\]
Similarly, Theorem~\ref{ad0-b} shows that 
\[
\#\AD_0\{\be\}=\{ x\in Q^\vee \mid (x,\ap)\in\{1,-1,-2,\ldots \} \ \ \forall \ap\in\Pi 
\quad \& \quad (x,\theta)\in \{0,2,3\ldots \} \}\ . 
\]
%%%%%%%%%%%%%%%
\end{comment}
%
By Theorems~\ref{ad-b} and \ref{ad0-b}, we have to count the points in $Q^\vee$
satisfying certain constraints.
%%In both cases, we exploit the equations of the corresponding simplex and
%%exclude the possibility that $(x,\ap)=0$ for $\ap\in\Pi$ and
%%$(x,\theta)=1$.
However, for practical computations, it is easier to deal with points of
the coweight lattice in $V$, denoted $P^\vee$. 
Let $\{\varpi_i\}_{i=1}^n$ be the basis for $V$ that is dual to
$\{\ap_i\}_{i=1}^n$. Then the lattice generated by the
$\varpi_i$'s is $P^\vee$.
If $y=\sum_{i} y_i\varpi_i\in P^\vee$, then $y\in Q^\vee$ if and only if
a certain congruence condition (depending on $\g$) is satisfied for $(y_1,\dots,y_n)\in
{\Bbb Z}^n$. 

Our primary goal is to compare the numbers \ 
$\#\{y\in D_{min}\cap Q^\vee \mid
   y\not\in H_\ap \ \ \forall\,\ap\in\HP \}$ and \ 
$\#\{y\in D_{min}\cap P^\vee \mid
   y\not\in H_\ap \ \ \forall\,\ap\in\HP \}$; and likewise for $D_{max}$.
Define the integers $c_i\in{\Bbb Z}$, $i=1,2,\ldots,n$, by the
formula $\theta=\sum_{i=1}^n c_i\ap_i$. 

It is clear that \ 
$\#\{y\in D_{min}\cap P^\vee \mid  y\not\in H_\ap \ \ \forall\,\ap\in\HP \}$
is equal to the number of solutions of the system of equations
\[
%%\begin{equation}  \label{syst}
    \left\{ \begin{array}{c}
          y_i\in\{-1,1,2,\ldots\} \ \ (i=1,2,\dots,n) \\
          c_1y_1+\ldots +c_n y_n\le 2   \\
           c_1y_1+\ldots +c_n y_n \ne 1
%%\\      \text{congruence determined by }\ \g  
\ \  .
\end{array}\right.
%%\end{equation}
\]
It is convenient to set $y_0=1-(c_1y_1+\ldots +c_n y_n)$. The new variable
$y_0$ also ranges over $\{-1,1,2,\ldots\}$, so that, letting $c_0=1$,
the above system takes a more symmetric form
\begin{equation}  \label{syst-ext}
    \left\{ \begin{array}{c}
          y_i\in\{-1,1,2,\ldots\} \ \ (i=0,1,\dots,n) \\
          c_0y_0+c_1y_1+\ldots +c_n y_n=1  
%%\\      \text{congruence determined by }\ \g  
\ \  .
\end{array}\right.
\end{equation}
In a sense, this procedure corresponds to taking the extended 
Dynkin diagram of $\g$. For this reason, system~\re{syst-ext} will
be referred to as the {\it min-extended\/} system.

Similarly, 
$\#\{y\in D_{max}\cap P^\vee \mid  y\not\in H_\ap \ \ \forall\,\ap\in\HP \}$
is equal to the number of solutions of the system of equations
\[
%%\begin{equation}  \label{syst0}
    \left\{ \begin{array}{c}
          y_i\in\{1,-1,-2,\ldots\} \ \ (i=1,2,\dots,n) \\
          c_1y_1+\ldots +c_n y_n\ge 0   \\
           c_1y_1+\ldots +c_n y_n \ne 1  
%%\\      \text{congruence determined by }\ \g  
\ \  .
\end{array}\right.
%%\end{equation}
\]
%%Replacing each $y_i$ with $-y_i$ and then \
Letting $y_0=1-(c_1y_1+\ldots +c_n y_n)$ as above, 
we obtain the {\it max-extended\/} system
\begin{equation}  \label{syst0-ext}
    \left\{ \begin{array}{c}
          y_i\in\{1,-1,-2,\ldots\} \ \ (i=0,1,\dots,n) \\
          c_0y_0+c_1y_1+\ldots +c_n y_n=1  
%%\\      \text{congruence determined by }\ \g  
\ \  .
\end{array}\right.
\end{equation}
Replacing each $y_i$ with $-y_i$ yields another system, which is sometimes 
more convenient to deal with:
%%Sometimes it will be convenient to replace the last system with 
%%the following system having the same number of solutions:
%
\begin{equation}  \label{syst00-ext}
    \left\{ \begin{array}{c}
          y_i\in\{-1,1,2,\ldots\} \ \ (i=0,1,\dots,n) \\
          c_0y_0+c_1y_1+\ldots +c_n y_n=-1  
%%\\      \text{congruence determined by }\ \g  
\ \  .
\end{array}\right.
\end{equation}
The following theorem shows that counting points in $P^\vee$ in place of
$Q^\vee$ does not lead us far away from our purpose.  The number
$f=[P^\vee: Q^\vee]$ is called the {\it index of connection\/} of $\Delta$.
It is also equal to $\#\{j\mid c_j=1\}$.

\begin{s}{Theorem}   \label{1/f}
\begin{itemize}
\item[\sf (i)] \ $\#\{x\in D_{min}\cap Q^\vee \mid
   x\not\in H_\ap \ \ \forall\,\ap\in\HP \}=\frac{1}{f}{\cdot}
\#\{x\in D_{min}\cap P^\vee \mid x\not\in H_\ap \ \ \forall\,\ap\in\HP \}$;
\item[\sf (ii)] \ $\#\{x\in D_{max}\cap Q^\vee \mid
   x\not\in H_\ap \ \ \forall\,\ap\in\HP \}=\frac{1}{f}{\cdot}
\#\{x\in D_{max}\cap P^\vee \mid x\not\in H_\ap \ \ \forall\,\ap\in\HP \}$.
\end{itemize}
\end{s}\begin{proof}
The argument amounts to a direct case-by-case verification.
For each simple Lie algebra $\g$, we look at the effect of the additional 
congruence condition imposed on systems \re{syst-ext} and \re{syst0-ext}.
One can define an action of the cyclic group ${\Bbb Z}_f$ on the set of solutions
of \re{syst-ext} and \re{syst0-ext} such that each ${\Bbb Z}_f$-orbit has cardinality
$f$ and contains a unique representative from $Q^\vee$.
Since the technical details are completely the same as in the proof of Theorem~5.5 in
\cite{mima}, we omit them.
\end{proof}%
Obviously, the number of solutions of \re{syst-ext} is equal to the coefficient
of $x$ in the Laurent series 
\[
  \prod_{i=0}^n (x^{-c_i}+x^{c_i}+x^{2c_i}+\ldots)=
\prod_{i=0}^n(\frac{x^{-c_i}}{1-x^{c_i}}-1) \ .
\]
We use the standard notation that $[x^a]F(x)$ denotes the coefficient of $x^a$ in the
Laurent series $F(x)$. Therefore, combining Theorems~\ref{ad-b} and
\ref{1/f}(i) we obtain
\begin{equation}  \label{equa-b}
     \#\AD\{\be\}=\frac{1}{f} [x]\prod_{i=0}^n(\frac{x^{-c_i}}{1-x^{c_i}}-1) \ .
\end{equation}
Similarly, starting with \re{syst00-ext} and  combining Theorems~\ref{ad0-b}
and \ref{1/f}(ii), we obtain
\begin{equation}  \label{equa00-b}
     \#\AD_0\{\be\}=\frac{1}{f} [x^{-1}]\prod_{i=0}^n(\frac{x^{-c_i}}{1-x^{c_i}}-1) \ .
\end{equation}
The Equations~\re{equa-b} and \re{equa00-b} show that the cardinalities in question 
depend only on the multiset of the coefficients of the highest root. Therefore, we can
already conclude that 
these cardinalities are equal for $\spn$ and $\sono$.

Now, we are ready to consider the case of classical Lie algebras.
For the future use, we introduce some notation for trinomial coefficients.
The coefficient $[x^0](x^{-1}+1+x)^n$ is called the {\it central trinomial\/},
denoted $\textsf{ct}_n$ and
$[x](x^{-1}+1+x)^n$ is called the {\it next-to-central trinomial\/},
denoted $\textsf{nct}_n$. In an explicit form, we have
\[
\textsf{ct}_n=\sum_{k\ge 0}\frac{ n!}{k! k! (n-2k)!} \quad \text{and} \quad
\textsf{nct}_n=\sum_{k\ge 0}\frac{ n!}{k! (k+1)! (n-2k-1)!} \ .
\]
{\bf 1)} \ $\g={\frak sl}_{n+1}$. \quad
Here all $c_i=1$ and $f=n+1$. Therefore 
\begin{multline}     \label{syst-sl}
\# \AD (\slno)\{ \be \}=\frac {1}{n+1} [x] (\frac {x^{-1}}{1-x}-1)^{n+1}= \\
\frac{1}{n+1}[x]\sum_{i=1}^{n+1} (-1)^{n+1-i}\genfrac{(}{)}{0pt}{}{n+1}{i}
\frac{x^{-i}}{(1-x)^i}=\\
\frac{1}{n+1}\sum_{i=1}^{n+1} (-1)^{n+1-i}\genfrac{(}{)}{0pt}{}{n+1}{i}
[x^{i+1}]\frac{1}{(1-x)^i}= \\
\frac{1}{n+1}\sum_{i=1}^{n+1} (-1)^{n+1-i}\genfrac{(}{)}{0pt}{}{n+1}{i}
\genfrac{(}{)}{0pt}{}{2i}{i+1}= \\
\sum_{i=1}^{n+1}(-1)^{n+1-i}\genfrac{(}{)}{0pt}{}{n}{i-1}C_i=
\sum_{j=0}^{n}(-1)^{n-j}\genfrac{(}{)}{0pt}{}{n}{j}C_{j+1} \ ,
\end{multline}
where $C_i:=\frac{1}{i+1}\genfrac{(}{)}{0pt}{}{2i}{i}$ is the $i$-th {\it Catalan
number}. 
It is well-known that the last expression gives $M_n$, the $n$-th {\it Motzkin
number}, see e.g. \cite[p.99]{CMR}.
There is a rich literature devoted to Motzkin numbers, where the
reader may find various definitions/interpretations of these numbers, see e.g.
\cite{martin},\,\cite{CMR},\,\cite{dsh},\,\cite[Ex.\,6.37,6.38]{rstan2}.
In \cite{mima}, it was shown that $M_n$ gives the number of
minimax ideals in $\AD(\slno)$.
\\[.6ex]
Similarly, 
\begin{multline}     \label{syst00-sl}
\# \AD_0(\slno)\{\be\}=\frac {1}{n+1} [x^{-1}] (\frac {x^{-1}}{1-x}-1)^{n+1}=... \\
%%\frac{1}{n+1}[x]\sum_{i=1}^{n+1} (-1)^{n+1-i}\genfrac{(}{)}{0pt}{}{n+1}{i}
%%\frac{x^{-i}}{(1-x)^i}=\\
\frac{1}{n+1}\sum_{i=1}^{n+1} (-1)^{n+1-i}\genfrac{(}{)}{0pt}{}{n+1}{i}
[x^{i-1}]\frac{1}{(1-x)^i}=... \\
%%\frac{1}{n+1}\sum_{i=1}^{n+1} (-1)^{n+1-i}\genfrac{(}{)}{0pt}{}{n+1}{i}
%%\genfrac{(}{)}{0pt}{}{2i}{i+1}= \\
\sum_{i=1}^{n+1}(-1)^{n+1-i}\genfrac{(}{)}{0pt}{}{n}{i-1}C_{i-1}=
\sum_{j=0}^{n}(-1)^{n-j}\genfrac{(}{)}{0pt}{}{n}{j}C_{j} \ .
\end{multline}
This time, the last expression yields $R_n$, the $n$-th {\it Riordan number},
see \cite[p.99]{CMR}.
Similarity of the expressions for $\AD(\slno)\{\be\}$ and $\AD_0(\slno)\{\be\}$ 
reflects the well-known relation $M_n=R_n+ R_{n+1}$.
We refer to \cite{CMR} for a discussion of Catalan, Motzkin, and Riordan
numbers.

{\bf 2)} \ $\g={\frak sp}_{2n}$ or ${\frak so}_{2n+1}$. \quad
Here $c_0=c_1=1$, $c_2=\ldots =c_n=2$, and $f=2$.  Therefore we have to
compute the coefficients

\begin{multline*}     %%\label{syst-sp}
\frac{1}{2}[x^{\pm 1}]\bigl(\frac {x^{-1}}{1-x}-1\bigr)^{2}
\bigl(\frac{x^{-2}}{1-x^2}-1\bigr)^{n-1}= 
\\
=\frac{1}{2}[x^{\pm 1}]\bigl(\frac{x^{-1}-1+x}{1-x}\bigr)^2{\cdot} 
\bigl(\frac{x^{-2}}{1-x^2}-1\bigr)^{n-1}= 
\\
=\frac{1}{2}[x^{\pm 1}]\frac{(x^{-1}+x^2)^2}{(1-x^2)^2}{\cdot}
\bigl(\frac {x^{-2}}{1-x^2}-1\bigr)^{n-1}=
\\
=\frac{1}{2}[x^{\pm 1}](x^{-2}+2x+x^4)\sum_{i=0}^{n-1}(-1)^{n-1-i}
\genfrac{(}{)}{0pt}{}{n-1}{i}\frac{x^{-2i}}{(1-x^2)^{i+2}} \ .
\end{multline*}
For the parity reasons for degrees, we need only the summand $2x$ in the first factor.
That is, the last expression equals
%%$\displaystyle 
\[
[x^{\pm 1}]\sum_{i=0}^{n-1}(-1)^{n-1-i}
\genfrac{(}{)}{0pt}{}{n-1}{i}\frac{x^{-2i+1}}{(1-x^2)^{i+2}} \ .
\]
For the coefficient of $x^{-1}$, we obtain 
\begin{multline*}
\#\AD_0(\sono)=\#\AD_0(\spn)=
\sum_{i=0}^{n-1}(-1)^{n-1-i}\genfrac{(}{)}{0pt}{}{n-1}{i}[x^{2i-2}]
\frac{1}{(1-x^2)^{i+2}}=\\
\sum_{i=1}^{n-1}(-1)^{n-1-i}\genfrac{(}{)}{0pt}{}{n-1}{i}[x^{i-1}]
\frac{1}{(1-x)^{i+2}}
=\sum_{i=1}^{n-1}(-1)^{n-1-i}\genfrac{(}{)}{0pt}{}{n-1}{i}
\genfrac{(}{)}{0pt}{}{2i}{i-1}=\\
(n-1)\sum_{i=1}^{n-1}(-1)^{n-1-i}\genfrac{(}{)}{0pt}{}{n-2}{i-1}C_i=
(n-1)M_{n-2}=\textsf{nct}_{n-1} \ .
\end{multline*}
The last equality is explained e.g. in \cite[Section\,5]{mima}.
\\[.6ex]
For the coefficient of $x$, we obtain 
\begin{multline*}
\#\AD(\sono)=\#\AD(\spn)=
\sum_{i=0}^{n-1}(-1)^{n-1-i}\genfrac{(}{)}{0pt}{}{n-1}{i}[x^{2i}]
\frac{1}{(1-x^2)^{i+2}}=\ldots \\
=\sum_{i=0}^{n-1}(-1)^{n-1-i}\genfrac{(}{)}{0pt}{}{n-1}{i}
\genfrac{(}{)}{0pt}{}{2i+1}{i} \ .
%%=\\
%%(n-1)\sum_{i=0}^{n-1}(-1)^{n-1-i}\genfrac{(}{)}{0pt}{}{n-1}{i-1}C_i=
%%(n-1)M_{n-2}
% \ .
\end{multline*}
Write temporarily $X_n$ for the last expression. The binomial transform 
of $\{X_n\}$ yields
\[
  \genfrac{(}{)}{0pt}{}{2n-1}{n-1}=\sum_{k=0}^{n-1}\genfrac{(}{)}{0pt}{}{n-1}{k}X_k \ .
\]
Comparing with Eq.~\re{new}, which is proved below, we see that $X_n=\textsf{dir}_n$, the
{\it number of directed animals of size\/} $n$.

Some other expressions for $\textsf{dir}_n$ are
\[
  \textsf{dir}_n=\textsf{ct}_{n-1}+\textsf{nct}_{n-1}=
%%[x](x^{-1}+1+x)^{n-1}+[x^0](x^{-1}+1+x)^{n-1}=
\sum_{q\ge 0}
\genfrac{(}{)}{0pt}{}{q}{[q/2]}\genfrac{(}{)}{0pt}{}{n-1}{q} \ .
\]
(See \cite[Eq.\,(27)]{anim} and \cite[5.16]{mima}.) It is also not hard to prove that
$\#\AD(\sono)\{\be\}-\#\AD_0(\sono)\{\be\}=\textsf{ct}_{n-1}$.
This leads to a simple expression of the Riordan numbers via trinomial coefficients:
\[
          R_n=\textsf{ct}_{n}-\textsf{nct}_{n} \ .
\]

{\bf 3)} \ $\g={\frak so}_{2n}$, $n\ge 4$. \quad
Here $c_0=c_1=c_2=c_3=1$, $c_4=\ldots =c_n=2$, and $f=4$.
Therefore we have to compute the coefficients

\begin{multline*}     
\frac{1}{4}[x^{\pm 1}]\bigl(\frac {x^{-1}}{1-x}-1\bigr)^{4}{\cdot}
\bigl(\frac{x^{-2}}{1-x^2}-1\bigr)^{n-3}= 
\\
=\frac{1}{4}[x^{\pm 1}]\frac{(x^{-1}+x^2)^4}{(1-x^2)^4}{\cdot}
\bigl(\frac{x^{-2}}{1-x^2}-1\bigr)^{n-3}=
\\
%%=\frac{1}{4}[x^{\pm 1}](x^{-1}+x^2)^2\frac{1}{(1-x^2)^2}(\frac {x^{-2}}{1-x^2}-1)^{n-1}=
%%\\
=\frac{1}{4}[x^{\pm 1}](x^{-4}+4x^{-1}+6x^2+4x^5+x^8)\sum_{i=0}^{n-3}(-1)^{n-3-i}
\genfrac{(}{)}{0pt}{}{n-3}{i}\frac{x^{-2i}}{(1-x^2)^{i+4}} \ .
\end{multline*}
For the parity reasons for degrees, we need only the summand $4x^5+4x^{-1}$ in the first factor.
That is, the last expression equals
\[
  [x^{\pm 1}]\sum_{i=0}^{n-3}(-1)^{n-3-i}
\genfrac{(}{)}{0pt}{}{n-3}{i}\frac{x^{-2i+5}+x^{-2i-1}}{(1-x^2)^{i+4}} \ .
\]
For the coefficient of $x^{-1}$ we obtain 
\[
  \#\AD_0(\sone)\{\be\}= 
\sum_{i=0}^{n-3}(-1)^{n-3-i}\genfrac{(}{)}{0pt}{}{n-3}{i}\Bigl[
\genfrac{(}{)}{0pt}{}{2i}{i-3}+\genfrac{(}{)}{0pt}{}{2i+3}{i}
\Bigr] \ .
\]
For the coefficient of $x$ we finally obtain 
\[
   \#\AD(\sone)\{\be\}= 
\sum_{i=0}^{n-3}(-1)^{n-3-i}\genfrac{(}{)}{0pt}{}{n-3}{i}\Bigl[
\genfrac{(}{)}{0pt}{}{2i+1}{i-2}+\genfrac{(}{)}{0pt}{}{2i+4}{i+1}
\Bigr] \ .
\]
As in the previous case, we have the relation
\[
  \#\AD(\sone)\{\be\}-\#\AD_0(\sone)\{\be\}=\textsf{ct}_{n-1} \ . 
\]
There is also a connection between the values for $\sone$ and $\sono$.
Namely, 
\[
   \#\AD(\sono)\{\be\}-\#\AD(\sone)\{\be\}=M_{n-2} \ .
\]
{\bf 4)} \ $\g=\GR{F}{4}$. Here a straightforward calculation shows 
$\#\AD\{\be\}=19$, $\#\AD_0\{\be\}=11$.
\\[.6ex]
 For $\g=\GR{E}{6}$, calculations based on Eq.~\re{equa-b}
and \re{equa00-b} show that
$\#\AD\{\be\}=111$ and $\#\AD_0\{\be\}=53$.

The case of $\GR{G}{2}$ is easy, and the cases of $\GR{E}{7}$
and $\GR{E}{8}$ are too difficult to do them by hand.

                %%%%
%%%%%%%%%%%%%%%%
%%%%%%%%%%%%%%%%   Section 5
%%%%%%%%%%%%%%%%
                %%%%

\section{The case of ${\g}={\frak sl}_{n+1}$}
\label{sl}
\setcounter{equation}{0}

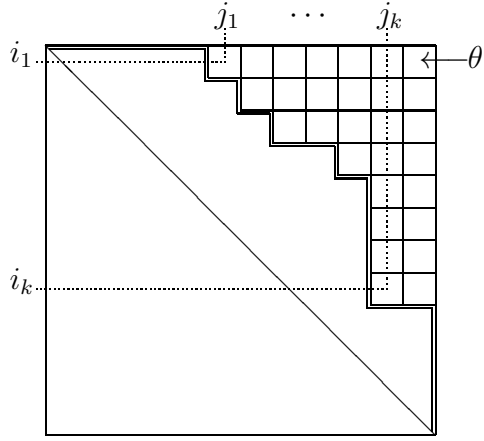
\begin{figure}[htb]
\setlength{\unitlength}{0.017in}
\begin{center}
\begin{picture}(125,125)(1,0)
\multiput(0,0)(120,0){2}{\line(0,1){120}}
\multiput(0,0)(0,120){2}{\line(1,0){120}}
\put(50,110){\line(1,0){70}}
\put(60,100){\line(1,0){60}}
\put(70,90){\line(1,0){50}}  % \put(95,105){\vector(0,-1){9}}
\put(90,80){\line(1,0){30}}  % \put(105,95){\vector(-1,0){9}}
\multiput(100,40)(0,10){4}{\line(1,0){20}}

\qbezier[39](-3,115),(30,115),(55,115)  % horizontal
\qbezier[76](-3,45),(50,45),(105,45)     % horizontal
\qbezier[50](105,125),(105,85),(105,45)
\qbezier[7](55,125),(55,120),(55,115)

%\put(50,70){\dashbox{1}(0,40){}}
%\put(50,70){\dashbox{1}(70,0){}}
%\put(40,80){\dashbox{1}(0,40){}}
%\put(40,80){\dashbox{1}(50,0){}}

\put(50,110){\line(0,1){10}}
\put(60,100){\line(0,1){20}}
\put(70,90){\line(0,1){30}}
\put(80,90){\line(0,1){30}}
\put(90,80){\line(0,1){40}}
\put(100,40){\line(0,1){80}}
\put(110,40){\line(0,1){80}}

\put(-11,115){$i_1$}
\put(-11,45){$i_k$}
\put(52,127){$j_1$}
\put(75,127){$\cdots$}
\put(102,127){$j_k$}
\put(130,113){$\theta$}
\put(115,113){$\longleftarrow$}

\put(120,0){\line(-1,1){120}}   %diagonal

\put(1,119){\line(1,0){48}}             \put(49,109){\line(0,1){10}}
\put(49,109){\line(1,0){10}}           \put(59,99){\line(0,1){10}}
\put(59,99){\line(1,0){10}}             \put(69,89){\line(0,1){10}}
\put(69,89){\line(1,0){20}}             \put(89,79){\line(0,1){10}}
\put(89,79){\line(1,0){10}}             \put(99,39){\line(0,1){40}}
\put(99,39){\line(1,0){20}}             \put(119,1){\line(0,1){38}}

\end{picture}
\end{center}
\caption{An {\sf ad}-nilpotent ideal in ${\frak sl}_{n+1}$}   \label{pikcha_A}
\end{figure}

\noindent
In this section, $\g={\frak sl}_{n+1}$.
%%and hence $p=n-1$. Here
We will explicitly describe $\AD\{\p\}$ for \un{ever}y
$\p\in {\frak Par}({\frak sl}_{n+1})$.
It will be shown that $\AD\{\p\}$ has a unique minimal element
and $\#\AD\{\p\}$ depends only on the difference $n-\srk\p$.
Using the duality construction from \cite[Section\,4]{duality}, we
produce a bijection between the minimax {\sf ad}-nilpotent ideals
and the ideals in $\AD\{\be\}$.
\\[.5ex]
We choose $\be$ (resp. $\te$) to be the space of upper-triangular (resp. diagonal)
matrices. With the usual numbering of rows and columns of matrices,
the positive roots are identified
with the pairs $(i,j)$, where $1\le i<j\le n+1$. For instance, $\ap_i=(i,i+1)$ and
therefore $(i,j)=\ap_i+\cdots +\ap_{j-1}$.
%%$\theta=(1,n)$.
An {\sf ad}-nilpotent ideal of $\be$ is represented by a right-justified
Ferrers diagram with at most $n$ rows, where the length of $i$-th row
is at most $n-i+1$. If a box of a Ferrers diagram corresponds to a
positive root $(i,j)$, then we say that this box has the coordinates
$(i,j)$.
The box containing the unique northeast corner of the diagram corresponds to $\theta$,
and the southwest corners give rise to the
generators of the corresponding ideal, see  Figure~\ref{pikcha_A}.

\noindent
An ideal $\ce$ (Ferrers diagram) is completely determined by the coordinates of boxes
%(roots)
that contain the southwest  corners of the diagram, say $(i_1,j_1),\dots,
(i_k,j_k)$. Then
we obviously  have
\[
\Gamma(\ce)=\{(i_1,j_1),\dots, (i_k,j_k)\}, \text{ where }\
1\le i_1< \dots  <i_k\le n,\ \
  2\le j_1< \dots  <j_k\le n{+}1 \ .
\]
Write $[n]$ for $\{1,2,\ldots,n\}$.
It is convenient to describe standard
parabolic subalgebras of $\slno$ by indicating the simple roots that
are not in the standard Levi subalgebra. That is, if
$E=\{l_1,l_2,\ldots,l_s\}$ is a subset of $[n]$,
then
\begin{equation}   \label{obozn-parab}
  \p(E)=\p(l_1,l_2,\ldots.l_s):=\p(\Pi\setminus\{\ap_{l_1},\ldots,\ap_{l_s}\}) \ .
\end{equation}
We always assume that $l_1<\ldots< l_s$.
Therefore the consecutive diagonal blocks of the standard Levi
subalgebra have sizes $l_1, l_2{-}l_1,\ldots, l_s{-}l_{s-1}, n+1-l_s$.
\\[.5ex]
Recall that $\#\AD(\slno)=\frac{1}{n+2}\genfrac{(}{)}{0pt}{}{2n+2}{n+1}=
C_{n+1}$, the $(n{+}1)$-th {\it Catalan number\/}, see e.g. \cite{CP1}.
There is a host of combinatorial objects counted by Catalan
numbers, see \cite[ch.\,6, Ex.\,6.19]{rstan2} and the ``Catalan addendum'' at
\verb|www-math.mit.edu/~rstan/ec|.
One of the main interpretations is that
$C_n$ is equal to
the number of lattice paths from $(0,0)$ to $(n,n)$ with steps $(1,0)$ and $(0,1)$,
%never rising above   the line $x=y$.
always staying in the domain $x\le y$, i.e., the number of Dyck paths of
semilength $n$. The Dyck path corresponding to a $\ce\in\AD(\slno)$
is the double path in Figure~\ref{pikcha_A}.

Recall that $M_s$ is the $s$-th {\it Motzkin number\/}. 
%%We refer to \cite{dsh} and \cite[Ex.\,6.37,6.38]{rstan2} for numerous
%%definitions/interpretations of these numbers.
One of the possible explicit expressions for them 
is
\begin{equation}   \label{motzkin}
M_s:=\sum_{r\ge 0}\genfrac{(}{)}{0pt}{}{s}{2r}C_r \ .
\end{equation}
\vskip-1ex

\begin{s}{Theorem}  \label{norm-p-sl}
\begin{itemize}
\item[\sf (1)] \
If\/ $\Gamma(\ce)=\{(i_1,j_1),\dots, (i_k,j_k)\}$,
then $\n_\g(\ce)=\p(\{i_1,\ldots,i_k\}\cup\{j_1{-}1,\ldots,j_k{-}1\})$;
\item[\sf (2)]
\ given a subset\/ $E=\{l_1,l_2,\ldots,l_s\}\subset [n]$,
the ideals in $\AD\{\p(E)\}$ are in a bijection
with the pairs $1\le a_1< \ldots <a_k\le n$, $1\le b_1< \ldots <b_k\le n$
of integer sequences such that $a_i\le b_i$ and
$\{a_1,\ldots,a_k\}\cup\{b_1,\ldots,b_k\}=E$.
(Here $k$ is not fixed, but it obviously satisfies the constraints
$s/2\le k\le s$.) The ideal corresponding to such a pair of sequences has the
generators $\{(a_i,b_i+1)\mid 1\le i\le k\}$.
\item[\sf (3)] \ If $\# E=s$, then
$\#\AD\{\p(E)\}=M_s$.
That is, this cardinality depends only on $s=\rk(\slno)-
\srk\p(E)$.
\end{itemize}
\end{s}\begin{proof*}
1. The simple roots that can be subtracted from the generators
have the numbers $i_1,\ldots,i_k,j_1{-}1,\ldots,j_k{-}1$ (repetitions are
allowed!). Then the description of $\n_\g(\ce)$ follows from Theorem~\ref{p-r}.

2. This readily follows from part (1).

3. Clearly, the number of pairs of sequences described in part (2)
depends only on $s$ and not on the explicit
form of $E$.
For instance, we may assume that $E=\{1,2,\ldots,s\}$. 
Then the assertion on the number of
the above pairs of sequences is precisely
the characterization of Motzkin numbers given in
\cite[Ex.\,6.38(e)]{rstan2}.
\par
Another (more "honest") way to see the connection with Motzkin numbers is
as follows. Let us temporarily write ${\goth M}_s$ for the cardinality
of $\#\AD\{\p(E)\}$. Since the number of standard parabolic
subalgebras of $\slno$ with semisimple corank $r$ equals $\genfrac{(}{)}{0pt}{}{n}{r}$,
$\AD(\slno)=\displaystyle\bigsqcup_\p \AD(\slno)\{\p\}$, and the cardinality
of $\AD(\slno)$ is known,
we obtain for each $n\in \Bbb N$ the identity
\[
       C_{n+1}=\sum_{r=0}^n\genfrac{(}{)}{0pt}{}{n}{r}{\goth M}_r \ .
\]
According to \cite{dsh}, there is an explicit relation between the Catalan and
Motzkin numbers which is exactly of such form. Hence ${\goth M}_r=M_r$
for all $r$.
\end{proof*}%
%
%
%There is an extensive literature devoted to Motzkin numbers.
%We only mention here papers \cite{martin},\,\cite{CMR},\,\cite{dsh}.
%In \cite{mima}, it was shown that $M_n$ gives the number of
%minimax ideals in $\AD(\slno)$.

\begin{s}{Theorem}  \label{min-and-max}
As above, $E=\{l_1,l_2,\ldots,l_s\}$.
The poset $\AD\{\p(E)\}$ has a unique maximal and unique
minimal {\sf ad}-nilpotent ideal. More precisely,
\begin{itemize}
\item \ $\ce_{max}(E)=\p(E)^{nil}$ and \\
$\Gamma(\ce_{max}(E))=\{(l_i, l_{i}+1) \mid 1\le i\le s\}$;
\item \ $\Gamma(\ce_{min}(E))=\{(l_i, l_{[s/2]+i}+1) \mid 
1\le i\le [(s{+}1)/2]\}$.
\end{itemize}
\end{s}\begin{proof*}
The first claim is a particular case of Lemma~\ref{fibre-psi}.
As for the second claim, it is clear that the ideal with given
generators lies in
$\AD\{\p(E)\}$.
Suppose $(a_1,\ldots,a_k)$, $(b_1, \ldots ,b_k)$ is an arbitrary pair
of sequences as in Theorem~\ref{norm-p-sl}(2).
Since each $l_j$ must appear among the $a_i$'s and $b_i$'s, we have $k\ge [(s+1)/2]$.
Also, $a_i\ge l_i$ and $b_i\le l_{s-k+i}\le l_{s-[(s+1)/2]+i}=
l_{[s/2]+i}$. These inequalities mean that each root $(a_i,b_i+1)$ lies in the
ideal with the generators $(l_j, l_{[s/2]+j}+1)$, $1\le j\le [(s{+}1)/2]$, 
which completes the proof.
%%The minimality follows by a direct verification.
\end{proof*}%
\begin{rem}{Remarks} 
Theorems~\ref{norm-p-sl}, \ref{min-and-max} have a number of interesting consequences.

1. Notice that $\ce_{min}(E)$ is always an Abelian ideal.
This reflects the fact that, for $\slno$ (and $\spn$), the mapping
$\ce\mapsto \n_\g(\ce)$ sets up a bijection between ${\frak Par}$ and the subset of
Abelian ideals in $\AD$, see \cite[Section 3]{pr}.
Thus, $\ce_{min}(E)$ is the unique Abelian ideal with normalizer
$\p(E)$.
The description of the minimal ideal in $\AD\{\p\}$ can also be stated in
a "coordinate-free" form. If $\srk\p=n-s$, then the minimal ideal in
$\AD\{\p\}$ is   
%%$\ce_{min}(\p)=(\p^{nil})^{[(s{+}1)/2]}$.
$(\p^{nil})^{[(s{+}1)/2]}$.

2. It is easily seen that the poset
$\AD(\slno)\{\p(l_1,l_2,\ldots,l_s)\}$
is naturally isomorphic to the poset $\AD({\frak sl}_{s+1})\{\be\}$.
In other words, the structure of $\AD\{\p(l_1,l_2,\ldots,l_s)\}$
depends only on the number of diagonal blocks of a Levi subalgebra, but not
on the sizes of blocks.
Therefore, in a sense, it suffices to consider only the {\sf ad}-nilpotent
ideals whose normalizer equals $\be$.

3. The decomposition $\AD(\slno)=\displaystyle\bigsqcup_\p \AD(\slno)\{\p\}$
yields a "materialization" of the identity
%\[
       $C_{n+1}=\sum_{r=0}^n\genfrac{(}{)}{0pt}{}{n}{r}M_r$.
%\]
\end{rem}%
As in \cite{mima}, we write $\AD_{mm}=\AD_{mm}(\g)$ for the subset of
minimax ideals in $\AD(\g)$. Recall that $\ce$ is called minimax, if $w_{min,\ce}=
w_{max,\ce}$, which means in particular that $\AD_{mm}\subset \AD_0$.
The geometric characterization is that $\ce$ is minimax if and only if
$R_\ce$ is a single alcove.
In \cite[Sect.\,6]{mima}, we obtained a description of the minimax ideals for
$\slno$ and $\spn$.  In particular, $\#\AD_{mm}(\slno)=M_n$.

Now, we show that the equality $\#\AD(\slno)\{\be\}=\#\AD_{mm}(\slno)$
is not a mere coincidence. To this end, recall the notion of the
{\it dual\/} ideal for an ideal $\ce\in\AD(\slno)$.
If $\Gamma(\ce)=\{(i_1,j_1),\dots, (i_k,j_k)\}$, where
$1\le i_1< \dots  <i_k\le n$ and
$2\le j_1< \dots  <j_k\le n+1$, then put \\[.5ex]
\centerline{
$X(\ce)=\{i_1,\ldots,i_k\}$ and $\tilde Y(\ce)=\{j_1{-}1,\ldots,j_k{-}1\}$.}
\\[.5ex]
By definition, the dual ideal for $\ce$, denoted $\ce^*$, is the ideal
determined by the equalities $X(\ce^*)=[n]\setminus \tilde Y(\ce)$
and $\tilde Y(\ce^*)=[n]\setminus X(\ce)$. The operation $\ce\mapsto\ce^*$
is well-defined, and $(\ce^*)^*=\ce$. It has also a number of other nice properties,
see \cite[Section\,4]{duality} for more details.

\begin{s}{Theorem}  \label{dual-sl}
For $\g=\slno$, we have
$\ce\in\AD_{mm}$ if and only if\/ $\n_\g(\ce^*)=\be$.
\end{s}\begin{proof}
Let $\ce\in\AD(\slno)$ be an arbitrary ideal with
$X(\ce)=\{a_1,\ldots,a_k\}$ and $\tilde Y(\ce)=\{b_1,\ldots,b_k\}$.
Then we have the following two characterizations:
\begin{itemize}
\item $\ce$ is minimax if and only if $X(\ce)\cap \tilde Y(\ce)=\varnothing$.
\ (See \cite[Corollary\,6.5]{mima});
\item $\n_\g(\ce)=\be$ if and only if $X(\ce)\cup \tilde Y(\ce)=[n]$.
\ (See Theorem~\ref{norm-p-sl}(2).)
\end{itemize}
The very definition of $\ce^*$ shows that
$X(\ce)\cap \tilde Y(\ce)=\varnothing$ if and only if
$X(\ce^*)\cup \tilde Y(\ce^*)=[n]$.
\end{proof}%
As a consequence, we immediately obtain

\begin{s}{Corollary}  \label{selfdual-sl }
An ideal $\ce$ is self-dual (i.e., $\ce=\ce^*$) if and only if
$\ce$ is minimax and $\n_\g(\ce)=\be$.
\end{s}%
%
%%{\bf Remark.} 
It was shown in \cite[Theorem\,4.6]{duality} that the number of 
self-dual ideals in $\AD({\frak sl}_{2m+1})$ is equal to $C_{m}$ 
(while it is obviously 0 for ${\frak sl}_{2m}$).
%%%\\[.7ex]
\begin{rem}{Remark}
Considering the normalizers of minimax ideals yields a materialization 
%%Let us obtain a proof ("materialization") 
of Eq.~\re{motzkin}.
%%using the normalizers of minimax ideals.
Since $\ce$ is minimax if and only if $X(\ce)\cap \tilde Y(\ce)=\varnothing$, we see
that $\ce\in\AD_{mm}$ has $k$ generators if and only if $\n_\g(\ce)$ has semisimple 
corank $2k$. 
(That is, unlike the general
case, there is a strong correlation between $\srk\n_\g(\ce)$
and $\#\Gamma(\ce)$.)
%%the number of generators of $\ce$.)
Obviously, any $\p\in {\frak Par}(\slno)$ having an even semisimple corank appears in 
this way. On the other hand,  if $\p=\p(l_1,\ldots, l_{2k})$, then
the minimax ideals in $\AD\{\p\}$ are in a bijection with 
the \un{dis}j\un{oint} partitions
\[
   \{l_1,\ldots,l_{2k}\}=\{a_1,\ldots,a_{k}\}\sqcup \{b_1,\ldots,b_{k}\}
\]
such that $a_1<\ldots< a_{k}$, $b_1<\ldots < b_{k}$, and $a_i < b_i$. 
It is well-known and easy to prove that the number of such partitions equals $C_k$.
Thus,
%%\\[.6ex]\centerline{
if $n-\srk\p=2k$, then $\#\{\ce\in\AD_{mm}\mid \n_\g(\ce)=\p \}=C_k$.
%%}
This yields Eq.~\re{motzkin}.
\end{rem}
\vskip-1.5ex

                %%%%
%%%%%%%%%%%%%%%%
%%%%%%%%%%%%%%%%   Section 6
%%%%%%%%%%%%%%%%
                %%%%

\section{The case of ${\g}={\frak sp}_{2n}$}
\label{sp}
\setcounter{equation}{0}

\noindent
Roughly speaking, the results for ${\frak sp}_{2n}$ are similar to those
for $\slno$. One of the notable distinctions is that 
the cardinalities of posets $\AD(\spn)\{\p\}$ are now expressed in terms of
{\it numbers of directed animals\/}
in place of Motzkin numbers.
\\[.5ex]
We use a standard matrix model of ${\frak sp}_{2n}$ corresponding
to a Witt basis for alternating bilinear form. 
For this basis of ${\Bbb C}^{2n}$, the algebra ${\frak sp}_{2n}$ has the 
following block form:
\[
{\frak sp}_{2n}=\{\left(\begin{array}{cr} A & B \\ C & D \end{array}
\right)\mid B={\widehat B},\ C={\widehat C},\ D=-{\widehat A} \} \ ,
\]
where $A,B,C,D$ are $n\times n$ matrices and $A\mapsto {\widehat A}$ is the 
transpose relative to the antidiagonal.
If $\ov{\be}$ is the standard Borel subalgebra of ${\frak sl}_{2n}$, then
$\be:=\ov{\be}\cap {\frak sp}_{2n}$ is a Borel subalgebra of ${\frak sp}_{2n}$.
(See also \cite[5.1]{duality}.) 
We identify the positive roots of ${\frak sp}_{2n}$ with 
the set $\{(i,j)\mid i<j,\ i+j\le 2n+1\}$.
Here the simple roots are $\ap_i=(i,i+1)$, $1\le i\le n$, and therefore:
\\[.5ex]\indent
\hphantom{$\bullet$ --} $(i,j)=\left\{\begin{array}{cc}
\ap_i+\ldots+\ap_{j-1}, & \text{if }\ j\le n+1 \ ,\\
\ap_i+\ldots +\ap_{2n-j}+2(\ap_{2n-j+1}+\ldots +\ap_{n-1})+\ap_n,
& \text{if }\ j > n+1 \ .\end{array}\right.$
\\[.7ex]
The root $(i,j)$ is long if and only if $i{+}j=2n+1$, i.e., the corresponding
matrix entry lies on the anti-diagonal.
The ideals for ${\frak sp}_{2n}$ can be identified with the ideals
(Ferrers diagrams) for ${\frak sl}_{2n}$ that are
symmetric with respect to the antidiagonal (=\,{\it self-conjugate}).
In other words, there is a natural bijection between
$\AD({\frak sp}_{2n})$ and the self-conjugate ideals in
$\AD({\frak sl}_{2n})$. More precisely,
suppose $\bar \ce\in \AD({\frak sl}_{2n})$ and
$\Gamma(\bar \ce)=\{(i_1,j_1),\ldots,(i_k,j_k)\}$
with $i_1<\ldots < i_k$, where we use our convention on the roots of
${\frak sl}_{2n}$. By definition, $\bar \ce$ is
self-conjugate if and only if $i_m+j_{k+1-m}=2n+1$ for all $m$.
Then the  corresponding ideal $\ce\in\AD({\frak sp}_{2n})$
has the generators $\Gamma(\ce)=\{(i_m,j_m)\mid m\le [(k+1)/2]\}$.
Conversely, given $\ce\in\AD(\spn)$, we obtain $\Gamma(\bar\ce)$ by replacing
each $(i,j)\in\Gamma(\ce)$ with $(i,j)$ and $(2n{+}1{-}j,2n{+}1{-}j)$.
If $\#\Gamma(\ce)=s$, then
\[\#\Gamma(\bar\ce)=
\left\{\begin{array}{ll} 2s, & \text{ if all roots in $\Gamma(\ce)$ are short}\\
                2s-1, & \text{ if $\Gamma(\ce)$ contains a (unique!) long root}.
\end{array}\right.
\]
%%Under this convention, any generator of an ideal in $\Delta^+({\frak sp}_{2n})$
%%satisfies the constraint $i_m+j_m\le 2n+1$.
We shall say that $\bar \ce\in\AD({\frak sl}_{2n})$ is the {\it symmetrization\/}
of $\ce\in \AD({\frak sp}_{2n})$. It is also clear that $\n_{\sltn}(\bar\ce)$
is a self-conjugate (in the above sense) standard
parabolic subalgebra, and that
$\n_{\sltn}(\bar\ce)$ is the symmetrization of $\n_{\spn}(\ce)$.
\\[.5ex]
Our general line of reasoning is as follows. Given an ideal $\ce\in\AD(\spn)$,
we take its symmetrization $\bar\ce\in\AD(\sltn)$
and then work with the generators of $\bar\ce$.
If $\Gamma(\bar \ce)=\{(i_1,j_1),\ldots,(i_k,j_k)\}$, then put
$E_\ce:=
(\{i_1,\ldots,i_k\}\cup\{j_1{-}1,\ldots,j_k{-}1\})\cap [n]
=(X(\bar\ce)\cup \tilde Y(\bar\ce))\cap [n]$.
We use the same notation for parabolic subalgebras as for $\slno$, see
Eq.~\re{obozn-parab}.

\begin{s}{Proposition}  \label{norm-p-sp}
\begin{itemize}
\item[\sf (i)] \
For any $\ce\in\AD(\spn)$, we have\/ $\n_\g(\ce)=\p(E_\ce)$.
\item[\sf (ii)] \ Given $E \subset [n]$, there is a bijection between the ideals in\/
$\AD\{\p(E)\}$ and the the pairs
$1\le a_1< \ldots <a_k\le 2n{-}1$, $1\le b_1< \ldots <b_k\le 2n{-}1$
of integer sequences such that $a_i\le b_i$, $a_i+b_{k+1-i}=2n$ for any  $i$, and
$(\{a_1,\ldots,a_k\}\cup\{b_1,\ldots,b_k\})\cap [n]=D$.
(Here $k$ is not fixed.)
The ideal corresponding to such a pair of sequences has the
generators $\{(a_i,b_i+1)\mid 1\le i\le [(k+1)/2]\}$.
\end{itemize}
\end{s}\begin{proof}
(i) By Theorem~\ref{p-r}, it suffices to realize which simple roots can be
subtracted from the generators of $\ce$.
For $(i,j)$ with $j\le n+1$ these are $\ap_i$ and $\ap_{j-1}$, while for
$(i,j)$ with $j > n+1$ these are $\ap_i$ and $\ap_{2n-j+1}$.
Thus, in view of previous definitions, the numbers of these roots form exactly
the set $E_\ce$.

(ii) This is just a reformulation of the above formulae expressing $E_\ce$ in terms
of the generators of the symmetrization of $\ce$.
\end{proof}%
Obviously, two strongly increasing sequences $(a_i)$, $(b_i)$, as in
Theorem~\ref{norm-p-sp}(ii), are determined by the intersections
$\{a_1,\ldots,a_k)\}\cap [n]$ and $\{b_1,\ldots,b_k\}\cap [n]$, cf. the next proof.
\\
Now, we show that the long simple root, $\ap_n$, plays a special role in
the symplectic case.

\begin{s}{Proposition}  \label{last-root}
Let $E$ be a subset of $[n{-}1]$. Then there is a natural bijection between
$\AD\{\p(E)\}$ and $\AD\{\p(E\cup\{n\})\}$.
\end{s}\begin{proof}
Suppose $\n_\g(\ce)=\p(E)$ and let $(a_1,\ldots,a_k)$, $(b_1,\ldots,b_k)$ be
the corresponding integer sequences as described in
Proposition~\ref{norm-p-sp}(ii). Then $n$ does not appear in both of them.
Let us insert $n$ in the appropriate place of each sequence.
It is easily seen that the two sequences obtained satisfy the
conditions of Proposition~\ref{norm-p-sp}(ii)
and thereby determine an ideal in $\AD\{\p(E\cup\{n\})\}$.
Indeed, we have
%%\begin{multline}
\begin{equation}   \label{2-sec}
\begin{array}{l}
a_1< \ldots < a_t< n< a_{t+1}< \ldots < a_k , \\
b_1< \ldots < b_{k-t}< n< b_{k-t+1}< \ldots < b_k ,
\end{array}
\end{equation}
%%\end{multline}
for some $t$
in view of the conditions $a_i+b_{k+1-i}=2n$. Furthermore,
$t\ge k-t$. For, otherwise we would have $a_{t+1}> n > b_{t+1}$.

Obviously, this procedure can be reversed.
With this notation, we also have
$E=\{a_1, \ldots, a_t\}\cup\{b_1, \ldots, b_{k-t}\}$, where the union is not
necessarily disjoint.
\end{proof}%
Thanks to this result, we can restrict ourselves to considering only
parabolic subalgebras $\p(E)$ with $E\subset [n{-}1]$.

\begin{s}{Theorem}  \label{zveri-bij}
Suppose $E=\{l_1,\ldots,l_s\}\subset [n{-}1]$. Then there is a bijection between
the ideals in $\AD\{\p(E)\}$ and the words $(v_1\ldots v_s)$
in the alphabet $\{-1,0,1\}$ such that all partial sums $\sum_{i\le m}v_i$
are non-negative.
\end{s}\begin{proof}
As we already know, an ideal in $\AD\{\p(E)\}$ is determined by two
sequences \re{2-sec} such that $a_i\le b_i$, $a_i+b_{k+1-i}=2n$,
and $\{a_1, \ldots, a_t\}\cup\{b_1, \ldots, b_{k-t}\}=E$.
The word $(v_1\ldots v_s)$ is defined by the following rule:
\\
\centerline{
$v_i=\left\{
\begin{array}{rl}
   1, & \text{if }\ l_i\in \{a_1, \ldots, a_t\}\setminus\{b_1, \ldots, b_{k-t}\} \\
   0, & \text{if }\ l_i\in \{a_1, \ldots, a_t\}\cap\{b_1, \ldots, b_{k-t}\} \\
  -1, & \text{if }\ l_i\in \{b_1, \ldots, b_{k-t}\}\setminus\{a_1, \ldots, b_{t}\} .
\end{array}\right.$
}
It is easily seen that non-negativity of all partial sums is guaranteed
by the condition $a_i\le b_i$.

Conversely, given a word $(v_1\ldots v_s)$, we can restore the presentation
of $E$ as a not necessarily disjoint
union $\{a_1, \ldots, a_t\}\cup\{b_1, \ldots, b_{k-t}\}$.
The non-negativity of partial sums implies that $t>k-t$ and $a_i\le b_i$.
Therefore these "truncated" $a$- and $b$-sequences
can be extended to the whole sequences \re{2-sec}, using
the conditions $a_i+b_{k+1-i}=2n$.
\end{proof}%
The number of {\it directed animals of size\/} $n$, denoted $\textsf{dir}_n$,
is defined as the number of certain $n$-element subsets of a two-dimensional
lattice. An explicit expression is
\begin{equation}  \label{dir-yavno}
  \textsf{dir}_n=\sum_{q\ge 0}\genfrac{(}{)}{0pt}{}{q}{[q/2]}\genfrac{(}{)}{0pt}{}{n-1}{q}
\ ,  
\end{equation}
see \cite[Eq.\,(27)]{anim}.
According to a  beautiful result of
D.~Gouyou-Beauchamps and G.~Viennot \cite[Theorem\,1]{anim},
the number of directed animals of size $s+1$
equals the number of words $(v_1\ldots v_s)$, as above.
Thus, we obtain the following

\begin{s}{Corollary}  \label{zveri-chislo}
If $E\subset [n-1]$ and $\#(E)=s$, then
\[
\#\AD(\spn)\{\p(E)\}=\#\AD(\spn)\{\p(E)\cup\{n\}\}=\textsf{dir}_{s+1} \ .
\]
In particular, $\#\AD(\spn)\{\be\}=\#\AD(\spn)\{\p(\ap_n)\}=\textsf{dir}_{n}$.
\end{s}%
Using this corollary and the equality $\#\AD(\spn)=
\genfrac{(}{)}{0pt}{}{2n}{n}$ \cite{CP1},
the decomposition $\AD(\spn)=\sqcup_\p \AD(\spn)\{\p\}$ is being translated to the
identity
\[
   \genfrac{(}{)}{0pt}{}{2n}{n}=2\sum_{k=0}^{n-1}\genfrac{(}{)}{0pt}{}{n{-}1}{k}
\textsf{dir}_{k+1}
\]
or
\begin{equation}  \label{new}
  \genfrac{(}{)}{0pt}{}{2n{-}1}{n-1}=
\sum_{k=0}^{n-1}\genfrac{(}{)}{0pt}{}{n{-}1}{k}\textsf{dir}_{k+1} \ ,
\end{equation}
which is seem to be new.
All explicit results on {\sf ad}-nilpotent ideals for $\spn$ are based
on the relation between an ideal $\ce\in\AD(\spn)$ and 
its symmetrization $\bar\ce\in\AD(\sltn)$.
There is an analogue of Proposition~\ref{min-and-max}:

\begin{s}{Proposition}  \label{mm-sp}
For any $\p\in {\frak Par}(\spn)$, the poset $\AD\{\p\}$ has a unique
minimal element (ideal), which is Abelian.
\end{s}\begin{proof}
Take the self-conjugate parabolic subalgebra
$\bar\p\in {\frak Par}(\sltn)$ corresponding to $\p$.
It follows from Theorem~\ref{min-and-max} that the minimal ideal in
$\AD({\frak sl}_{2n})\{\bar\p\}$ is self-conjugate and
Abelian, and therefore it determines a (necessarily minimal and
Abelian) ideal in $\AD(\spn)\{\p\}$.
\end{proof}%
Recall the construction of the dual ideal in the symplectic setting.
Given $\ce\in\AD(\spn)$, one defines $\ce^*$ in a round-about way via
the symmetrization. That is, $\ce^*$ is the ideal whose
symmetrization is $(\bar\ce)^*\in\AD(\sltn)$, see \cite[Sect.\,5]{duality}.

\begin{s}{Theorem}  \label{dual-sp}
For $\g=\spn$, we have
$\ce\in\AD_{mm}$ if and only if\/ $\n_\g(\ce^*)=\be$.
\end{s}\begin{proof}
It follows from Proposition~\ref{norm-p-sp}(ii) that
$\n_\g(\ce)=\be$ if and only if $\n_{\sltn}(\bar\ce)=\bar\be$.
On the other hand, the description of the minimax ideals in $\AD(\spn)$
\cite[Cor.\,6.8]{mima} essentially says that $\ce\in\AD_{mm}(\spn)$ if and
only if $\bar\ce\in\AD_{mm}(\sltn)$. Thus, the result follows from
Theorem~\ref{dual-sl}.
\end{proof}%
Finally, we consider the normalizers of minimax ideals and characterize the sets
$\AD\{\p\}\cap\AD_{mm}$.

\begin{s}{Theorem}  \label{norm-mm-sp}
\begin{itemize}
\item[\sf (i)] \ The algebra $\p(E)\in {\frak Par}(\spn)$ is the normalizer of a minimax
ideal if and only if $E\subset [n{-}1]$.
\item[\sf (ii)] \ For $E\subset [n{-}1]$ with $\# E=s$, we have \
%%\centerline{
$\#\{\ce\in\AD_{mm} \mid \n_\g(\ce)=\p(E)\}=\genfrac{(}{)}{0pt}{}{s}{[s/2]}$.
\end{itemize}
\end{s}\begin{proof}
(i) If $\bar\ce\in\AD(\sltn)$ is the symmetrization of a minimax ideal $\ce$, then
the condition $X(\bar\ce)\cap\tilde Y(\bar\ce)=\varnothing$, which characterizes
the minimax ideals, readily implies that $n\not\in X(\bar\ce)\cup\tilde Y(\bar\ce)$.

(ii) If $E=\{l_1,\ldots,l_s\}$, then the minimax ideals in $\AD(\p(E))$ are in
a bijection with the \un{dis}j\un{oint} partitions
\[
     \{a_1,\ldots,a_t\}\sqcup \{b_1,\ldots,b_{s-t}\}=E
\]
such that $a_1<\ldots <a_t$, $b_1<\ldots <b_{s-t}$, and $a_i \le b_i$.
Using the rule from the proof of Theorem~\ref{zveri-bij}, one
sees that the number of such partitions is equal to the number of words
$(v_1\ldots v_s)$ in the alphabet $\{-1,1\}$ such that all partial sums
$\sum_{i\le m}v_i$ are nonnegative. A direct calculation shows that the latter is equal
to $\genfrac{(}{)}{0pt}{}{s}{[s/2]}$, see Lemma~\ref{s/2}.
\end{proof}%
Notice that this theorem and the partition
$\AD_{mm}(\spn)=\sqcup_\p \AD_{mm}(\spn)\{\p\}$ yield a materialization of the 
identity Eq.~\re{dir-yavno}. 
\\
For convenience of the reader, we give a proof of the following result, which was used
in the proof of Theorem~\ref{norm-mm-sp}. 

\begin{s}{Lemma}   \label{s/2}
The number of words $(v_1\ldots v_s)$ in the alphabet $\{-1,1\}$ 
such that all partial sums
$\sum_{i\le m}v_i$ are nonnegative is equal to $\genfrac{(}{)}{0pt}{}{s}{[s/2]}$.
\end{s}\begin{proof}
Consider all words $(v_1\ldots v_s)$ such that the total sum 
$\sum_{i\le s}v_i$ is nonnegative. Clearly, the number of such words is
\[
  \genfrac{(}{)}{0pt}{}{s}{[s/2]}+\genfrac{(}{)}{0pt}{}{s}{[s/2]{-}1}
  +\genfrac{(}{)}{0pt}{}{s}{[s/2]{-}2}+\ldots \ .
\]
(The number of $-1$'s can be any integer $\le [s/2]$.)
Let us say that a word is {\it bad\/} if at least one partial sum is negative.
%%Our aim is to count the number of bad words. 
Suppose $(v_1\ldots v_s)$ is bad,
and let $\sum_{i\le t}v_i=-1$ be the \un{first} negative partial sum.
Consider the word $(w_1\ldots w_s)$, where
$w_i=\left\{\begin{array}{rc} -v_i, & \text{if } \ i\le t , \\
                               v_i, & \text{if } \ i> t  .      
\end{array}\right.$
Then $\sum_{i\le s}w_i=\sum_{i\le s}v_i+2\ge 2$, and it is easily seen that this
procedure yields a bijection between the bad words and all the words with the total
sum $\ge 2$. The number of the latter is
\[
  \genfrac{(}{)}{0pt}{}{s}{[s/2]{-}1}+\genfrac{(}{)}{0pt}{}{s}{[s/2]{-}2}
  +\ldots \ .
\]
Whence the number of non-bad words is equal to $\genfrac{(}{)}{0pt}{}{s}{[s/2]}$.
\end{proof}%
The reader may recognize that we have used in the proof
the reflection principle for lattice paths in an algebraic form.

                %%%%
%%%%%%%%%%%%%%%%
%%%%%%%%%%%%%%%%   Section 7
%%%%%%%%%%%%%%%%
                %%%%

\section{Some remarks on $\slno$, $\spn$, and other simple Lie algebras}
\label{other}
\setcounter{equation}{0}

\noindent
It seems that $\slno$ and $\spn$ are the most attractive simple Lie
algebras from the point of view of the theory of {\sf ad}-nilpotent ideals. 
Let us list some relevant nice properties of these two series
that do not hold in general:

\begin{itemize}
\item \ There is a natural procedure of constructing the dual {\sf ad}-nilpotent
ideal, see \cite[Sect.\, 4 \& 5]{duality}.
\item \ Taking the dual ideal 
yields a bijection between $\AD_{mm}$ and $\AD\{\be\}$,
see Theorems~\ref{dual-sl} and \ref{dual-sp}.
\item \ There is an explicit description of the minimax ideals,
see \cite[Sect.\,6]{mima}.
\item \ Write $\Ab=\Ab(\g)$ for the set of all Abelian ideals.
It is shown in \cite{pr} that the assignment $\ah\mapsto \n_\g(\ah)$,
$\ah\in\Ab$, yields a bijection between $\Ab(\g)$ and ${\frak Par}(\g)$
if and only if $\g=\sln$ or $\spn$.
\end{itemize}
As we have seen in Sections~\ref{sl} and \ref{sp}, 
considering various classes of ideals for 
$\slno$ and $\spn$  provides a tool for demonstrating some identities
related to Catalan and Motzkin numbers and numbers of directed animals.
Let us present one more speculation of this kind.
Consider the generating function for the number of minimax ideals in $\g$
with a given semisimple corank of the normalizer:
\[
{\mathcal F}_{n\text{-}mm}(\g,t):=\sum_{s} \#\{\ce\in\AD_{mm} \mid
\rk\g-\srk\n_\g(\ce)=s\}{\cdot}t^s \ ,
\]
Our computation in Theorem~\ref{norm-mm-sp} shows that
\[
{\mathcal F}_{n\text{-}mm}(\spn,t)=\sum_{s=0}^{n-1} \genfrac{(}{)}{0pt}{}{n-1}{s}
\genfrac{(}{)}{0pt}{}{s}{[s/2]}t^s \ ,
\]
and we know that ${\mathcal F}_{n\text{-}mm}(\spn,1)=\textsf{dir}_n$.
It is curious to observe that ${\mathcal F}_{n\text{-}mm}(\spn,-1)=R_n$,
the $n$-th Riordan number. 
%%(See \cite{CMR} about these numbers.)
\\[.6ex]
In case of $\sono$, there is an ersatz construction of duality that is based
on the similarity of shifted  Ferrers diagrams representing {\sf ad}-nilpotent
ideals in $\spn$ and $\sono$ \cite[Sect.\,5]{duality}.
However, this construction does not yield a bijection 
between $\AD_{mm}(\sono)$ and $\AD(\sono)\{\be\}$, which is already seen for
$n=3$. Also, no explicit
description of minimax ideals for $\sono$ is known.

Furthermore, the cardinalities of the sets $\AD_{mm}$ and $\AD\{\be\}$ are not always
equal. The following table, which presents results of our explicit calculations,
shows various possibilities:
\begin{center}
\begin{tabular}{c|ccccc|}
 &  ${\frak so}_8$ &  ${\frak so}_{10}$ & $\GR{E}{6}$ & $\GR{F}{4}$ & $\GR{G}{2}$ 
\\ \hline
$\#\AD_{mm}$     &  9  &  23 &  67 & 17 & 3 \\
$\#\AD\{\be\}$   &  11 &  31 & 111 & 19 & 2 
%%$\#\AD_0\{\be\}$ &  4  &  12 &  53 & 11 & 1 
\end{tabular}
\end{center}

\end{document}